\documentclass[12pt,francais]{article} 
\usepackage{amsmath,amsthm,amssymb,amscd}
\usepackage{graphicx}
\usepackage{babel}
\evensidemargin = \oddsidemargin
\advance\textwidth by -2in
\def\ly{\fontencoding{U}\fontfamily{lasy}\fontseries{m}\fontshape{n}\selectfont}
\def\guil#1{\leavevmode\hbox{{\ly(\kern-0.20em(\kern+0.20em}}\nobreak{}\,#1\,%
  \nobreak\hbox{{\ly\kern+0.20em)\kern-0.20em)}}}
\theoremstyle{plain}
\newtheorem{theoreme}{Th\'eor\`eme}[section] 
\newtheorem{proposition}[theoreme]{Proposition} 
\newtheorem{lemme}[theoreme]{Lemme}
\newtheorem{corollaire}[theoreme]{Corollaire} 
\newtheorem*{assertion}{Assertion}
\theoremstyle{definition}
\newtheorem{definition}[theoreme]{D\'efinition} 
\theoremstyle{remark}
\newtheorem*{remarque}{Remarque} 
\newtheorem*{remarques}{Remarques}
\newtheorem*{exemple}{Exemple}

\def\Alinea#1{\hfill\break%
  \hbox to \parindent{\hss{\textup{#1}}\enspace}\ignorespaces}
\def\alinea#1{\noindent%
  \hbox to \parindent{\hss{\textup{#1}}\enspace}\ignorespaces}
\def\legende#1{\par \vskip10pt \parbox{\textwidth}{#1}}
\def\up{\textup}
\def\from{\colon}
\def\res#1{\,\vert\,{}_{#1}}
\def\abs#1{\lvert #1 \rvert}
\def\mod#1{\mathrel{(\mathrm{mod}}#1)}
\def\classe#1{\mathcal{C}_{}^{#1}}
\def\Chi{\setbox0=\hbox{$\chi$} \mathord{\raise\dp0\hbox{$\chi$}}}
\def\eps{\varepsilon}
\def\C{\mathbf{C}}
\def\D{\mathbf{D}}
\def\H{\mathbf{H}}

\def\L{\mathbf{L}}
\def\N{\mathbf{N}}
\def\P{\mathbf{P}}

\def\R{\mathbf{R}}
\def\S{\mathbf{S}}
\def\T{\mathbf{T}}
\def\Z{\mathbf{Z}}
\def\RP^#1{\mathbf{P}^{#1}(\mathbf{R})}
\def\HH{\mathcal{H}}
\def\id{\mathrm{id}}
\def\PSL{\mathrm{PSL}}
\def\Card{\mathop{\mathrm{Card}}}
\def\Im{\mathop{\mathrm{Im}}\nolimits}
\def\Adh{\mathop{\mathrm{Adh}}\nolimits}
\def\Int{\mathop{\mathrm{Int}}\nolimits}
\def\Aut{\mathop{\mathrm{Aut}}\nolimits}
\def\pgcd{\mathop{\mathrm{pgcd}}\nolimits}
\def\dist{\mathop{\mathrm{dist}}\nolimits}
\def\aire{\mathop{\mathrm{aire}}\nolimits}
\def\l{\mathtt{l}}
\def\tb{\mathtt{tb}}
\def\q{\mathtt{i}}
\def\w{\mathtt{e}}
\def\op{\mathopen}
\def\cl{\mathclose}
\def\ol{\bar}

\def\wt{\tilde}
\def\wh{\hat}
\def\adresse#1{\def\Y{\egroup\egroup\hbox\bgroup\itshape\bgroup} \par \noindent
  \hbox to \textwidth {\hfill
  \vbox {\small \hbox \bgroup\itshape\bgroup #1 \egroup\egroup}}}
\title{Structures de contact sur les vari\'et\'es\\
fibr\'ees en cercles au-dessus d'une surface}
\author{Emmanuel \textsc{Giroux}%
\thanks{\ Centre National de la Recherche Scientifique (\textsc{umr} 5669)}}
\date{Octobre 1999}
\begin{document}
\maketitle

Dans cet article, on essaie d'analyser le comportement global des structures de
contact sur les vari\'et\'es fibr\'ees en cercles au-dessus d'une surface close.
Plusieurs \'etudes ant\'erieures motivent et guident ce travail. Tout d'abord,
sur les fibr\'es principaux en cercles, les structures de contact invariantes
admettent, \`a isotopie \'equivariante pr\`es, une classification remarquable
\cite{Lu}. D'autre part, de nombreux r\'esultats sur les repr\'esentations du
groupe fondamental d'une surface dans $\PSL_2(\R)$ et les hom\'eomorphismes du
cercle \cite{EHN,Gh:euler,Gh:fuchs,Ma,Mi,Wo} contribuent \`a mettre \`a jour la
structure topologique et dynamique des feuilletages de codimension $1$ sur les
vari\'et\'es de dimension~$3$ fibr\'ees en cercles (voir par exemple \cite{Le,
Th:godbillon,Th:these}). Or, en dimension~$3$, les structures de contact sont,
avec les fibr\'es tangents des feuilletages de codimension~$1$, les seuls champs
de plans localement homog\`enes et les d\'eveloppements parall\`eles des deux
th\'eories ont fait appara\^\i tre de nombreux traits communs.

\smallskip

Soit $V$ une vari\'et\'e connexe et orient\'ee, fibr\'ee en cercles au-dessus
d'une surface~$S$ close, orientable et de caract\'eristique d'Euler $\Chi(S)$
n\'egative ou nulle. Avant de pr\'esenter rapidement les principaux r\'esultats
de ce travail, on rappelle que~$V$, en tant que vari\'et\'e lisse orient\'ee,
est identifi\'ee par le nombre d'Euler $\Chi(V,S)$ de la fibration $\pi \from V
\to S$ \ (cf. section 1.B). On rappelle aussi qu'une structure de contact
(directe) sur~$V$ est un champ de plans d\'efini localement comme le noyau d'une
$1$-forme~$\alpha$ dont le produit ext\'erieur avec $d\alpha$ est une forme
volume positive pour l'orientation choisie.

Dans la partie~1, on d\'emontre que $V$ admet une \guil{connexion strictement
convexe}, c'est-\`a-dire une structure de contact (directe) transversale aux
fibres, si et seulement si le nombre d'Euler $\Chi(V,S)$ est inf\'erieur ou
\'egal \`a~$-\Chi(S)$ \ (th\'eor\`eme~\ref{t:milnorwood}). Cette in\'egalit\'e
est un reliquat de l'in\'egalit\'e de Milnor-Wood~\cite{Mi,Wo}, laquelle peut en
retour \^etre interpr\'et\'ee comme suit : il existe sur~$V$ une connexion plate
--~\emph{i.e.} un feuilletage transversal aux fibres~-- si et seulement si
cohabitent sur~$V$ des structures de contact directes et indirectes
transversales aux fibres.
 
Dans la partie~2, on donne une caract\'erisation g\'eom\'etrique des structures 
de contact qui sont isotopes \`a des connexions. Pr\'ecis\'ement, on prouve que 
toute structure de contact orientable~$\xi$ sur~$V$ satisfait l'alternative
exclusive suivante (th\'eor\`eme~\ref{t:thurston}) : ou bien $\xi$ est isotope
\`a une connexion, ou bien il existe, dans un rev\^etement fini de $(V,\xi)$,
une courbe legendrienne isotope \`a la fibre et le long de laquelle $\xi$ ne
tourne pas, \emph{i.e.} d\'etermine la m\^eme trivialisation normale que la
fibration~$\pi$. Ce r\'esultat est la version de contact d'un th\'eor\`eme
d\'emontr\'e par W.~{\sc Thurston} dans~\cite{Th:these} (voir aussi \cite{Le}),
selon lequel un feuilletage sur~$V$ est isotope \`a une connexion (plate) si et
seulement si la fibre n'est pas isotope \`a une courbe trac\'ee sur une feuille.

Dans la partie~3, on classifie les structures de contact transversales aux
fibres \`a isotopie et conjugaison pr\`es (th\'eor\`eme~\ref{t:ghys}). On montre
en particulier que deux connexions strictement convexes peuvent ne pas \^etre
isotopes. La formule suivante r\'esume bien la situation : il y a deux sortes de
structures de contact transversales aux fibres, celles qui sont tangentes aux
fibres et les autres. Ces derni\`eres existent d\`es que $\Chi(V,S) < -\Chi(S)$
et appartiennent \`a une m\^eme classe d'isotopie. Les premi\`eres, en revanche,
n'existent que si $n \Chi(V,S) = -\Chi(S)$ pour un certain entier $n\ge1$ et
leur classification se ram\`ene \`a celle des rev\^etements fibr\'es \`a $n$
feuillets de $V$ au-dessus de la vari\'et\'e $\S(TS)$ des droites orient\'ees
tangentes \`a~$S$ ; on montre ainsi que, pour $\Chi(S) < 0$, elles forment
autant de classes de conjugaison qu'il y a de diviseurs de~$n$. 

Dans la partie~4, on tente de classer les structures de contact qui ne sont pas
isotopes \`a des connexions. Comme le sort des structures de contact vrill\'ees
est scell\'e par un th\'eor\`eme de Y.~{\sc Eliashberg}~\cite{El:overtwist}, on
s'int\'eresse aux structures de contact tendues. Celles-ci pr\'esentent des
comportements tr\`es diff\'erents selon qu'elles sont virtuellement vrill\'ees ou
universellement tendues (cf. d\'efinition~\ref{d:eliashberg}). Les structures de
contact virtuellement vrill\'ees sur~$V$ constituent un nombre fini de classes
d'isotopie (th\'eor\`eme~\ref{t:bennequin}) born\'e par $1 + \sup \{0,-\Chi(S)
-\Chi(V,S)-1\}$. En revanche, les structures de contact universellement tendues
et non isotopes \`a des connexions forment une infinit\'e de classes d'isotopie
qui sont en bijection naturelle avec les classes d'isotopie de syst\`emes (non
vides) de courbes essentielles sur~$S$ (th\'eor\`eme~\ref{t:lutz}). En fait,
ces structures sont toutes isotopes \`a des structures de contact invariantes
(par une quelconque action libre du cercle qui d\'efinit la fibration) et leur
classification \`a isotopie pr\`es co\"{\i}ncide avec la classification des
structures invariantes \`a isotopie \'equivariante pr\`es.

En parall\`ele avec ce dernier r\'esultat, on \'etablit une in\'egalit\'e de
Bennequin semi-locale (proposition~\ref{p:bennequin}) qui conduit \`a la
classification des structures de contact tendues et $\R$-invariantes sur le
produit par $\R$ de toute surface~$F$ close et orientable (th\'eor\`eme \ref
{t:convexe}) : les classes d'isotopie de ces structures sont \`a nouveau en
bijection avec les classes d'isotopie de syst\`emes de courbes essentielles sur
$F$. L'int\'er\^et de ce r\'esultat tient au fait qu'une surface $F$ plong\'ee
dans une vari\'et\'e de contact de dimension $3$ poss\`ede g\'en\'eriquement un
voisinage tubulaire trivialis\'e $U \cong F \times \R$ dans lequel la structure
de contact est $\R$-invariante \cite{Gi:convexite}.

\smallskip

Je tiens \`a remercier ici \'Etienne {\sc Ghys}, Jean-Pierre {\sc Otal} et Bruno
{\sc S\'evennec} avec qui j'ai eu de nombreuses discussions sur certains aspects
de ce travail. D'autre part, Fran\c cois {\sc Lalonde} et Dietmar~{\sc Salamon}
m'ont offert l'occasion de pr\'esenter les r\'esultats discut\'es dans ce texte
\`a Montr\'eal en juin~1995 et \`a Warwick en mars~1998 ; je les en remercie
vivement.

\section{Existence de structures de contact transversales}

\subsection{Comment contacter Milnor-Wood}

Soit~$V$ une vari\'et\'e orient\'ee, fibr\'ee en cercles au-dessus d'une surface
close~$S$. Le th\'eor\`eme ci-dessous relate ce qui reste de l'in\'egalit\'e de
Milnor-Wood~\cite{Mi,Wo} quand on cherche non pas des feuilletages transversaux
aux fibres --~connexions plates~-- mais des structures de contact (directes)
transversales aux fibres --~connexions \guil{strictement convexes}.

\begin{theoreme} \label{t:milnorwood}
Soit $V$ une vari\'et\'e connexe orient\'ee, fibr\'ee en cercles au-dessus d'une
surface close~$S$. Pour que $V$~porte une structure de contact transversale aux
fibres, il faut et il suffit que le nombre d'Euler $\Chi(V,S)$ de la fibration
$V \to S$ v\'erifie l'in\'egalit\'e
$$ \left\{ \begin{aligned}
\Chi(V,S) & \le -\Chi(S) &\quad& \text{si \ $\Chi(S) \le 0$,} \\ 
\Chi(V,S) & < 0 &\quad& \text{si \ $\Chi(S) > 0$.}
\end{aligned} \right. $$
\end{theoreme}

\smallskip

\begin{remarques}
\Alinea{a)}
La d\'efinition du nombre d'Euler $\Chi(V,S)$, et en particulier de son signe en
fonction de l'orientation de~$V$, est rappel\'ee dans la section~B. Pour une
vari\'et\'e~$V$ fibr\'ee en cercles orient\'es au-dessus d'une surface~$S$
orient\'ee, $\Chi(V,S)$ co\"{\i}ncide avec le nombre de Chern du fibr\'e en
droites complexes associ\'e et l'orientation de~$V$ choisie est la juxtaposition
des orientations de la base et de la fibre. Ainsi, la vari\'et\'e des droites
orient\'ees tangentes \`a~$S$ a pour nombre d'Euler $-\Chi(S)$ lorsqu'elle est
munie de l'orientation induite par sa structure de contact canonique.
 
\alinea{b)}
Le th\'eor\`eme~\ref{t:milnorwood} a \'et\'e ind\'ependamment obtenu (du moins
pour une surface~$S$ orientable) par A.~{\sc Sato} et T.~{\sc Tsuboi}~\cite{ST}.
La preuve qu'on donne ci-apr\`es consiste simplement \`a adapter les arguments
de J.~{\sc Wood}~\cite{Wo}. Par ailleurs, vu comme espace total d'un fibr\'e
\emph{principal} en cercles, $V$ admet des connexions invariantes qui sont des
structures de contact directes si et seulement si $\Chi(V,S)$ est strictement
n\'egatif~\cite{Lu}.
\end{remarques}

\subsection{Nombre d'Euler et connexions}

On rappelle d'abord ce qu'est le nombre d'Euler $\Chi(V,S)$ de la fibration $\pi
\from V \to S$. On trace sur~$S$ un bouquet~$K$ de $2-\Chi(S)$ cercles ayant
pour compl\'ementaire un disque et on note~$D$ (resp.~$W$) le disque polygonal
(resp. le tore plein poly\'edral orient\'e) qu'on obtient en d\'ecoupant~$S$
(resp.~$V$) le long de~$K$ (resp. $\pi^{-1}(K)$). Au-dessus de~$K$, la fibration
$\pi$ admet des sections et chacune d'elles d\'etermine, sur le bord orient\'e
de~$W$, une courbe~$C$ \emph{dont la classe d'isotopie est invariable}. De
m\^eme, le bord des disques m\'eridiens de~$W$ est une courbe~$B$ bien
d\'efinie \`a isotopie pr\`es. Le nombre d'Euler $\Chi(V,S)$ est l'intersection
homologique $B \cdot C$ de~$B$ et~$C$, ces deux courbes \'etant orient\'ees de
mani\`ere \`a couper les fibres dans le m\^eme sens.

\smallskip

Si~$S$ est une surface orient\'ee, de genre~$g$, et si~$V$ est munie d'une
connexion~$\xi$ (champ de plans transversal aux fibres), le nombre d'Euler
$\Chi(V,S)$ s'interpr\`ete comme suit.

L'holonomie de~$\xi$ associe \`a chaque cercle orient\'e~$K_i$ du bouquet~$K$,
$1 \le i \le 2g$, un diff\'eomorphisme~$\phi_i$ de la fibre~$\S^1$ qui surplombe
le sommet de~$K$ : c'est l'application de premier retour qu'on obtient en
suivant les courbes int\'egrales de~$\xi$ au-dessus de~$K_i$. De
plus, chaque classe d'homotopie de sections de $\pi \res {K_i}$ d\'etermine un 
rel\`evement~$\wt\phi_i$ de~$\phi_i$ \`a~$\R$ : elle trivialise en effet $\pi
\res {K_i}$ si bien que le segment de courbe int\'egrale qui joint $1 \in \S^1$ 
\`a $\phi_i(1)$ se projette en un chemin sur la fibre~$\S^1$ ; on prend alors 
pour $\wt\phi_i(0)$ l'extr\'emit\'e dans~$\R$ du relev\'e partant de~$0$. 

De m\^eme, l'holonomie de la connexion induite par~$\xi$ sur~$W$, encore not\'ee
$\xi$, associe au bord orient\'e de~$D$ un diff\'eomorphisme du cercle qui, \`a
conjugaison pr\`es, s'\'ecrit comme un mot
$$
\phi = w \bigl( \phi_1^{}, \phi_1^{-1}, \dots,
   \phi_{2g}^{}, \phi_{2g}^{-1} \bigr) $$
dans lequel chaque~$\phi_i$, tout comme son inverse, intervient une fois et une
seule --~le d\'ecoupage de~$S$ d\'edouble chaque cercle orient\'e~$K_i$ en deux
ar\^etes de~$\partial D$ ayant des orientations incompatibles. Les courbes $B$
et~$C$, respectivement fournies par les sections de~$\pi$ au-dessus de $D$ et de
$K$, d\'eterminent alors deux rel\`evements distincts $\wt\phi_D$ et~$\wt\phi_K$
de~$\phi$ \`a~$\R$. Compte tenu de ce qui pr\'ec\`ede, ceux-ci v\'erifient les
identit\'es suivantes :
\begin{itemize}
\item
$\wt\phi_K = w \bigl( \wt\phi_1^{}, \wt\phi_1^{-1}, \dots, \wt\phi_{2g}^{}, \wt
\phi_{2g}^{-1} \bigr)$ o\`u les $\wt\phi_i$ sont des rel\`evements quelconques
des~$\phi_i$ ;
\item
$\wt\phi_D(t) - \wt\phi_K(t) = \Chi(V,S)$ pour tout r\'eel~$t$.
\end{itemize}

\subsection{In\'egalit\'es cl\'es}

\begin{proposition} \label{p:wood}
Soit $\xi$ une connexion sur~$V$ et $\wt\phi_K$, $\wt\phi_D$ les
diff\'eomorphismes de $\R$ d\'efinis plus haut.

\alinea{\textbf{a)}}
Quelle que soit la connexion~$\xi$,
$$ -2g \le \wt\phi_K(t) - t \le 2g \quad \text{pour tout r\'eel $t$.} $$
 
\alinea{\textbf{b)}}
Si~$\xi$ est une structure de contact directe,
$$ \wt\phi_D(t) - t < 0 \quad \text{pour tout r\'eel $t$.} $$
\end{proposition}

\begin{proof}
\Alinea{a)}
C'est l'in\'egalit\'e de J.~{\sc Wood}~\cite{Wo}. Soit~$\wt\HH$ l'espace des
hom\'eomorphismes de~$\R$ qui commutent avec la translation $t \mapsto t+1$.
L'application $h \from \wt\HH \to \R$ d\'efinie par
$$ h (\wt\psi) = \sup \bigl\{ \wt\psi(t) - t, \ t \in \R \bigr\} $$
v\'erifie
$$ h (\wt\psi_1 \wt\psi_2) \le h (\wt\psi_1) + h (\wt\psi_2)
   \le h (\wt\psi_1 \wt\psi_2) + 1 \quad
   \text{pour tous \ } \wt\psi_1, \wt\psi_2 \in \wt\HH \,. $$
On en d\'eduit la majoration
$$ h (\wt\phi_K)
 \le \sum_{i=1}^{2g} \Bigl( h (\wt\phi_i^{}) + h (\wt\phi_i^{-1}) \Bigr)
 \le 2g $$
et on obtient la minoration en observant que
$$ \inf \bigl\{ \wt\psi(t) - t, \ t \in \R \bigr\}
 = - \sup \bigl\{ \wt\psi^{-1}(t) - t, \ t \in \R \bigr\} \,. $$

\smallskip

\alinea{b)}
Il s'agit de feuilleter~$W = D \times \S^1$ par des disques m\'eridiens dont le
bord orient\'e soit partout transversal \`a~$\xi$ et pointe du m\^eme c\^ot\'e
que les fibres orient\'ees. Par approximation, on peut supposer que $D$~est un
carr\'e $[0,1]^2$. On choisit une coordonn\'ee~$\theta$ sur la fibre orient\'ee
en~$(0,0)$ et on l'\'etend \`a~$W$ en la d\'ecr\'etant constante sur chaque
courbe int\'egrale de~$\xi$ qui rev\^et soit un segment vertical soit la base
du carr\'e. Dans les coordonn\'ees $(x,y,\theta) \in [0,1]^2 \times \S^1$, la
connexion~$\xi$ a pour \'equation
$$ d\theta - u(x,y,\theta)\, dx = 0 \quad \text{o\`u \ 
   $u(x,0,\theta) = 0$ \ pour tout \ $(x,\theta) \in [0,1] \times \S^1$.} $$
La condition qui exprime alors que~$\xi$ est une structure de contact directe
par rapport \`a $dx \wedge dy \wedge d\theta$ s'\'ecrit $\partial_y u < 0$. Les
niveaux de~$\theta$ sont ainsi des m\'eridiens dont le bord orient\'e a la 
transversalit\'e souhait\'ee le long du c\^ot\'e $y=1$ --~car $u(x,1,\theta)<0$
pour tout $(x,\theta)$~-- et est tangent \`a~$\xi$ ailleurs. On peut donc les
perturber comme voulu.
\end{proof}

\medskip

Pour une surface~$S$ orientable et une connexion~$\xi$ qui est une structure de
contact directe, la proposition ci-dessus \'etablit l'in\'egalit\'e 
$$ \Chi(V,S) < 2g = 2-\Chi(S), $$
qui n'est le r\'esultat d\'esir\'e que si~$S$ est la sph\`ere. Si $S$~est le
plan projectif, le passage au rev\^etement double permet aussi de conclure. Dans
les autres cas, l'in\'egalit\'e du th\'eor\`eme~\ref{t:milnorwood} s'obtient
via l'astuce classique suivante. On prend un rev\^etement \`a $n$~feuillets
de~$S$ par une surface connexe orientable~$S_n$ et on note $V_n \to S_n$ le
rappel du fibr\'e $V \to S$ au-dessus de~$S_n$. Les relations
$$ \Chi(S_n) = n \Chi(S), \quad
   \Chi(V_n,S_n) = n \Chi(V,S) \quad \text{et} \quad
   \Chi(V_n,S_n) < 2 - \Chi(S_n) $$
donnent
$$ \Chi(V,S) < \frac2n - \Chi(S) \,.$$
Comme on peut choisir~$n$ arbitrairement grand, on obtient
$$ \Chi(V,S) \le -\Chi(S) \,. $$

\smallskip

\begin{remarque}
Si la connexion~$\xi$ est \emph{plate} au sens o\`u elle s'int\`egre en un
feuilletage, le diff\'eomorphisme~$\wt\phi_D$ est l'identit\'e et les arguments 
qui pr\'ec\`edent d\'emontrent l'in\'egalit\'e classique de Milnor-Wood, \`a
savoir
$$ \bigl\lvert \Chi(V,S) \bigr\rvert \le \sup \bigl\{ 0, -\Chi(S) \bigr\} \,. $$
Suite au travail de S.~{\sc Altschuler}~\cite{Al}, W.~{\sc Thurston} a invent\'e
le terme de \emph{feuilletact} (\emph{foliatact} dans~\cite{Th:feuilletacts}
mut\'e en \emph{confoliation} dans~\cite{ET}) pour d\'esigner un champ de plans
dont toute \'equation de Pfaff~$\alpha$ est telle que la $3$-forme $\alpha
\wedge d\alpha$ ne change pas de signe. La preuve de la proposition~\ref
{p:wood}-b montre que, si la connexion~$\xi$ est un
feuilletact direct ($\alpha \wedge d\alpha \ge 0$), le diff\'eomorphisme~$\wt
\phi_D$ v\'erifie 
$$ \wt\phi_D(t) - t \le 0 \quad \text{pour tout r\'eel $t$,} $$
de sorte que le nombre d'Euler $\Chi(V,S)$ satisfait \`a l'in\'egalit\'e
$$ \Chi(V,S) \le \sup \bigl\{ 0, -\Chi(S) \bigr\} \,. $$
\end{remarque}

\subsection{Construction de structures transversales}

\begin{lemme} \label{l:chirurgie}
Soit $V$ et~$V'$ deux vari\'et\'es connexes et orient\'ees, fibr\'ees en cercles
au-dessus d'une surface close $S$. Si $\Chi(V',S) \le \Chi(V,S)$ et si~$V$ porte
une structure de contact directe et transversale aux fibres, alors $V'$ en admet
une aussi.
\end{lemme}

\begin{proof}
On suppose pour simplifier que la surface~$S$ et les fibres sont orient\'ees.
La d\'efinition du nombre d'Euler donn\'ee dans la section~B montre que $V'$
s'obtient \`a partir de~$V$ par la chirurgie suivante : on retire \`a $V$ la
pr\'eimage $W \simeq \D^2 \times \S^1$ d'un disque de~$S$ et on recolle un autre
tore plein~$W'$ par un diff\'eomorphisme $\partial W' \to \partial W$ qui
respecte les fibres orient\'ees et envoie le bord de chaque disque m\'eridien
de~$W'$ sur une courbe de type $\bigl(1, \Chi(V,S) - \Chi(V',S) \bigr)$ dans le
produit $\partial W = \partial\D^2 \times \S^1$.

Si~$V$ porte une structure de contact~$\xi$ directe et transversale aux fibres,
la proposition~\ref{p:wood}-b montre qu'on peut feuilleter~$W$ par des disques
m\'eridiens dont le bord orient\'e soit transversal \`a~$\xi$ et pointe du
m\^eme c\^ot\'e que les fibres orient\'ees. Si $\Chi(V',S) \le \Chi(V,S)$, un
tel feuilletage par disques m\'eridiens existe aussi sur~$W'$ et il est alors
facile de prolonger \`a~$W'$ la structure $\xi \res{V \setminus W}$ pour obtenir
sur~$V'$ une structure de contact directe et transversale aux fibres.
\end{proof}

\medskip

Pour chaque surface~$S$, il reste donc \`a construire une structure de contact
directe et transversale aux fibres sur la vari\'et\'e orient\'ee~$V$ dont le
nombre d'Euler $\Chi(V,S)$ est le plus grand autoris\'e par l'in\'egalit\'e du
th\'eor\`eme~\ref{t:milnorwood}.

Pour la sph\`ere~$\S^2$, la structure de contact usuelle sur~$\S^3$ est
orthogonale aux cercles de Hopf et fournit l'exemple voulu. Pour toutes les
autres surfaces (y compris $\P^2$), la vari\'et\'e~$V$ est celle des droites
orient\'ees tangentes \`a~$S$. Elle porte une structure de contact canonique
$\xi_S$ qui induit l'orientation pour laquelle $\Chi(V,S) = -\Chi(S)$ mais qui
est tangente aux fibres et non pas transversale. En effet, si $\delta$ est une
droite orient\'ee tangente \`a~$S$ en un point~$q$, le plan $\xi_S (q,\delta)
\subset T_{(q,\delta)}V$ est l'image inverse de~$\delta$ par la projection. En
choisissant dans $\xi_S(q,\delta)$, pour tout $(q,\delta)$, un vecteur non
vertical dont la projection sur~$S$ donne l'orientation de $\delta$, on fabrique
sur~$V$ un champ de vecteurs legendrien non singulier et transversal aux fibres.
La condition de contact assure que, si on pousse~$\xi_S$ par le flot de ce champ
pendant un bref instant, on obtient une structure transversale aux fibres. Cet
argument d\'emontre plus g\'en\'eralement le fait suivant :

\begin{proposition} \label{p:perturbation}
Soit $V$ une vari\'et\'e connexe et orient\'ee, fibr\'ee en cercles au-dessus
d'une surface orientable~$S$. Toute structure de contact tangente aux fibres et
orientable le long des fibres est d\'eformable en une structure de contact
transversale aux fibres par une isotopie arbitrairement petite.
\end{proposition}

\subsection{Autre construction par la g\'eom\'etrie hyperbolique}

Pour clore cette partie, voici une autre construction adapt\'ee de la th\'eorie
des feuilletages. Soit $S$ une surface close orientable, de genre $g\ge1$, et
$V$ une vari\'et\'e connexe orient\'ee fibr\'ee en cercles au-dessus de~$S$.
D'apr\`es la discussion des sections B et~C, construire sur~$V$ une structure de
contact transversale aux fibres revient \`a trouver $2g$~diff\'eomorphismes
$\phi_i$ du cercle dont les rel\`evements $\wt\phi_i$ \`a $\R$ v\'erifient
$$ \bigl( \prod_{i=1}^g [\wt\phi_{2i-1}, \wt\phi_{2i}] \bigr) (t) - t <
 - \Chi(V,S) \quad \text{pour tout r\'eel $t$.} $$
L'argument ci-dessous, souffl\'e par \'E.~{\sc Ghys} et tout empreint de \cite
{Th:godbillon} (voir aussi \cite{EHN}), fournit $2g$~\'el\'ements $\phi_i$ de
$\PSL_2(\R)$ dont les rel\`evements~$\wt\phi_i$ v\'erifient
$$ \bigl( \prod_{i=1}^g [\wt\phi_{2i-1}, \wt\phi_{2i}] \bigr) (t) - t < \Chi(S)
   \quad \text{pour tout r\'eel $t$.} $$
Soit $P$ un polygone convexe \`a $4g$~c\^ot\'es dans le plan hyperbolique~$\H^2
$. On suppose que les sommets de~$P$, num\'erot\'es $s_1, \dots, s_{4g}$ dans le
sens des aiguilles d'une montre, v\'erifient
$$ \left\{ \begin{aligned}
   \dist (s_{4i-3},s_{4i-2}) &= \dist (s_{4i-1},s_{4i}) \\
   \dist (s_{4i-2},s_{4i-1}) &= \dist (s_{4i},s_{4i+1})
\end{aligned} \right. \quad
   \text{pour $1 \le i \le g$ et $s_{4g+1} = s_1$.} $$
On colle alors isom\'etriquement chaque ar\^ete orient\'ee $[s_{4i-3},s_{4i-2}]$
(resp. $[s_{4i-2},s_{4i-1}]$) sur l'ar\^ete orient\'ee $[s_{4i},s_{4i-1}]$
(resp.  $[s_{4i+1},s_{4i}]$). On obtient ainsi une surface close orientable~$S$
de genre~$g$ munie d'une m\'etrique hyperbolique ayant une singularit\'e conique
en~$s$, point image des sommets de~$P$. On pose ensuite $R = S \setminus \{s\}$,
on choisit dans $R$ un point de r\'ef\'erence~$r$, image d'un point $r_*\in \Int
P$ situ\'e tr\`es pr\`es de~$s_1$, et on note $(\wt R, \wt r)$ le rev\^etement
universel de $(R,r)$. Ces donn\'ees d\'eterminent une application d\'eveloppante
$D \from (\wt R, \wt r) \to (\H^2,r_*)$ et une repr\'esentation d'holonomie $h
\from \pi_1(R) \to \PSL_2(\R)$.

Soit $C \subset S$ le cercle trigonom\'etrique de centre~$s$ passant par~$r$
et $\gamma \in \pi_1(S)$ sa classe d'homotopie. Le point $D(\gamma \cdot \wt r)
= h(\gamma)\, (r_*)$ est l'image de $r_*$ par la rotation hyperbolique de centre
$s_1$ et d'angle la somme des angles int\'erieurs de~$P$, \`a savoir $(4g-2) \pi
- \aire(P)$ d'apr\`es la formule de Gauss-Bonnet. D'autre part, $h(\gamma)$ est
le produit de $g$~commutateurs dans $\PSL_2(\R)$. Pour les identifier, on note
que l'image inverse de~$C$ dans~$P$ est form\'ee de $4g$ arcs de cercles qui, en
partant de~$r_*$, sont centr\'es successivement aux points $s_1, s_4, s_3, s_2,
s_5, s_8, s_7, s_6, \dotsc$. Soit alors $\phi_{2i-1}$ et $\phi_{2i}$, $1 \le i
\le g$, les \'el\'ements de $\PSL_2(\R)$ caract\'eris\'es par les propri\'et\'es
suivantes :
\begin{align*}
\phi_{2i-1} (s_{4i-1}) &= s_{4i-2} &\quad \phi_{2i} (s_{4i-2}) &= s_{4i+1} \\
\phi_{2i-1} (s_{4i}) &= s_{4i-3}   &\quad \phi_{2i} (s_{4i-3}) &= s_{4i} .
\end{align*}
Par construction, $\prod_{i=1}^g [\phi_{2i-1},\phi_{2i}]$ vaut bien $h(\gamma)$.
Reste \`a d\'eterminer le nombre de translation du produit des commutateurs des
rel\`evements $\wt\phi_i$. Pour cela, on regarde le cas limite o\`u $P$ est un
polygone euclidien dans le plan tangent \`a $\H^2$ en un point~$s_0$. Dans ce
cas, les transformations~$\phi_i$ sont toutes des rotations de centre~$s_0$ et
commutent donc. Par suite, $\prod_{i=1}^g [\wt\phi_{2i-1},\wt\phi_i] = \id$. Il
en r\'esulte que, dans le cas g\'en\'eral, le nombre de translation vaut $-\frac
1{2\pi} \aire(P)$ et prend ainsi n'importe quelle valeur entre $0$ et $(1-2g)$.
En particulier, pour obtenir une structure de contact sur le fibr\'e $\S(TS)$
des droites orient\'ees tangentes \`a~$S$, il faut partir d'un polygone d'aire
sup\'erieure \`a $(4g-4)\pi$, \emph{i.e.} d'une m\'etrique ayant un atome de
courbure positive en la singularit\'e conique.

\section{Caract\'erisation des structures de contact trans\-versales}

\subsection{Comment contacter Thurston}

$V$~d\'esigne toujours une vari\'et\'e connexe orient\'ee, fibr\'ee en cercles
au-dessus d'une surface close~$S$. Dans~\cite{Th:these}, W.~{\sc Thurston} met
en \'evidence l'alternative exclusive suivante : si $\xi$ est un feuilletage de
codimension~$1$ orientable sur~$V$, ou bien $\xi$ est isotope \`a un feuilletage
transversal aux fibres, ou bien il existe une courbe simple trac\'ee sur une
feuille qui est isotope \`a la fibre. Dans le second cas, $\xi$ poss\`ede en
fait un ensemble minimal vertical --~\`a isotopie pr\`es~-- qui, pour peu que
$\xi$ soit $\classe2$ et que $S$ ne soit pas un tore, est n\'ecessairement une
feuille torique. Le r\'esultat qui suit est un analogue de ce th\'eor\`eme pour
les structures de contact. Son \'enonc\'e requiert un peu de terminologie.

\begin{definition} \label{d:eliashberg}
Soit~$\xi$ une structure de contact sur une vari\'et\'e~$M$ de dimension~$3$. On
dit que~$\xi$ est \emph{vrill\'ee} s'il existe un disque~$D$ plong\'e dans~$M$
qui est tangent \`a~$\xi$ en tous les points de son bord --~le disque~$D$ est
lui-m\^eme appel\'e \emph{disque vrill\'e}. On dit que $\xi$ est \emph{tendue}
si elle n'est pas vrill\'ee.

Une structure vrill\'ee sur~$M$ induit une structure encore vrill\'ee sur tout
rev\^etement de~$M$ mais une structure tendue peut aussi induire une structure
vrill\'ee sur certains rev\^etements. On dira qu'une structure de contact~$\xi$
sur~$M$ est \emph{virtuellement vrill\'ee} (resp. \emph{universellement tendue})
si elle est tendue et si elle induit une structure vrill\'ee (resp. tendue) sur
un rev\^etement fini de~$M$ (resp. sur le rev\^etement universel de~$M$). Ces
deux propri\'et\'es s'excluent mutuellement mais il n'est pas clair qu'elles
soient exactement compl\'ementaires l'une de l'autre. C'est cependant le cas si
le groupe fondamental de~$M$ est r\'esiduellement fini, donc par exemple
pour~$V$.
\end{definition}

Le long d'une fibre de~$V$, les champs de vecteurs normaux dont la projection
sur~$S$ est constante d\'eterminent une classe d'homotopie canonique de sections
non singuli\`eres du fibr\'e normal. De plus, d\`es que $V$ diff\`ere de $\S^2
\times \S^1$, cette classe d'homotopie est invariante par tout diff\'eomorphisme
de~$V$ isotope \`a l'identit\'e qui pr\'eserve la fibre consid\'er\'ee. Du coup,
le fibr\'e normal de toute courbe ferm\'ee simple isotope \`a la fibre poss\`ede
aussi une classe d'homotopie canonique de sections partout non nulles. Par abus
de langage, lorsque $V \simeq \S^2 \times \S^1$, une courbe isotope \`a la fibre
d\'esigne dans la suite une courbe munie d'une isotopie qui l'am\`ene sur une
fibre.

\begin{definition} \label{d:bennequin}
Soit $\xi$ une structure de contact sur~$V$. Pour toute courbe legendrienne~$L$
isotope \`a la fibre, on appelle \emph{enroulement de~$\xi$ autour de~$L$} --
ou \emph{enroulement de~$L$} -- le nombre $\w(L)$ de tours que $\xi$~fait le
long de~$L$ par rapport au champ normal canonique. On appelle \emph{enroulement}
de~$\xi$ et on note $\w(\xi)$ le supremum des enroulements $\w(L)$ pour toutes
les courbes legendriennes~$L$ isotopes aux fibres. Ce nombre est un entier si et
seulement si la structure de contact~$\xi$ est orientable le long des fibres.
\end{definition}

\smallskip

\begin{remarque}
D\`es lors que $V$ n'est ni un espace lenticulaire ni un tore, un th\'eor\`eme
de F.~{\sc Waldhausen}~\cite{Wa} assure que tous les diff\'eomorphismes de~$V$
respectent la fibration $\pi \from V \to S$, \`a isotopie pr\`es. L'enroulement
des structures de contact sur~$V$ est alors invariant par conjugaison.
\end{remarque}

\smallskip

Enfin, l'\'enonc\'e ci-dessous tient tacitement compte du fait que tout
rev\^etement fini de~$V$ fibre naturellement en cercles au-dessus d'un
rev\^etement fini de~$S$.

\begin{theoreme} \label{t:thurston}
Soit $V$ une vari\'et\'e connexe et orient\'ee, fibr\'ee en cercles au-dessus
d'une surface close et orientable~$S$. Toute structure de contact~$\xi$ sur~$V$
qui est orientable le long des fibres v\'erifie l'alternative exclusive
suivante\up{ :}
\begin{itemize}
\item
ou bien $\xi$ est isotope \`a une structure de contact transversale aux
fibres\up{ ;}
\item
ou bien il existe, dans un rev\^etement fini $(\wt V,\wt\xi)$ de $(V,\xi)$, une 
courbe legendrienne isotope \`a la fibre de~$\wt V$ et d'enroulement nul. Mieux,
le passage \`a un rev\^etement fini n'est n\'ecessaire que lorsque $\xi$ est
virtuellement vrill\'ee.
\end{itemize}
\end{theoreme}

Ce th\'eor\`eme est une cons\'equence des deux propositions ci-dessous dont la
d\'emons\-tration occupe la suite de la partie~~2.

\begin{proposition} \label{p:calculs}
Soit $\xi$ une structure de contact sur~$V$ orientable le long des fibres.

\alinea{\textbf{a)}}
L'enroulement $\w(\xi)$  appartient \`a $\Z \cup \{+\infty\}$. En outre, pour
tout entier $n \le \w(\xi)$, il existe dans $V$ une courbe legendrienne~$L$
isotope \`a la fibre et d'enroulement $\w(L)$ \'egal \`a~$n$.

\alinea{\textbf{b)}}
Si $\xi$ est vrill\'ee, son enroulement $\w(\xi)$ est infini.

\alinea{\textbf{c)}}
Si $\xi$ est transversale aux fibres, $\xi$ est universellement tendue et son
enroulement $\w(\xi)$ est strictement n\'egatif.
\end{proposition}

\begin{proposition} \label{p:thurston}
Soit $\xi$ une structure de contact sur~$V$ orientable le long des fibres. Si
$\xi$ est universellement tendue et si son enroulement $\w(\xi)$ est strictement
n\'egatif, $\xi$ est isotope \`a une structure de contact transversale aux
fibres.
\end{proposition}

\smallskip

\begin{proof}[D\'emonstration du th\'eor\`eme~\ref{t:thurston}]
Soit $\xi$ une structure de contact sur~$V$ orientable le long des fibres. La
proposition~\ref{p:calculs} montre tout d'abord que l'alternative envisag\'ee
pour $\xi$ est exclusive. Elle assure aussi que, si $\xi$ est vrill\'ee (resp.
virtuellement vrill\'ee), il existe dans~$V$ (resp. dans un rev\^etement fini
de~$V$) une courbe legendrienne isotope \`a la fibre et d'enroulement nul. Si
$\xi$ est au contraire universellement tendue, de deux choses l'une : ou bien
son enroulement $\w(\xi)$ est strictement n\'egatif et $\xi$ est isotope \`a une
structure transversale aux fibres (proposition~\ref{p:thurston}), ou bien
$\w(\xi)$ est positif ou nul et il existe dans~$V$ une courbe legendrienne
isotope \`a la fibre et d'enroulement nul (proposition~\ref{p:calculs}).
\end{proof}

\subsection{Estimations d'enroulement}

On d\'emontre ici la proposition~\ref{p:calculs}. Pour cela, on rappelle que
l'\emph{invariant de Thurston-Bennequin} $\tb(L)$ d'une courbe legendrienne~$L$
homologiquement nulle dans $(V,\xi)$ est l'enlacement de $L$ avec $L+\nu$, o\`u
$\nu$ est un champ de vecteurs normal \`a~$\xi$ le long de~$L$. L'enroulement
est un cousin de cet invariant ; en particulier, si $\Chi(V,S) = \pm1$, toute
courbe legendrienne~$L$ isotope \`a la fibre est homologiquement nulle et $\w(L)
= \tb(L) \pm 1$.

\smallskip

\alinea{a)}
Comme $\xi$ est orientable le long des fibres, elle effectue un nombre entier de
tours autour de chaque courbe legendrienne isotope \`a la fibre\footnote
{On observe au passage qu'une structure de contact transversale aux fibres est
automatiquement (co)\,orientable le long des fibres.}.
Par suite, $\w(\xi) \in \Z \cup \{+\infty\}$. Soit mainte\-nant $L_0$ une courbe
legendrienne isotope \`a la fibre et d'enroulement $\w(L_0) \ge n$, $n \in \Z$.
Soit d'autre part $L_1$ un n\oe ud legendrien topologiquement trivial contenu
dans une boule~dis\-jointe de $L_0$ et dont l'invariant de Thurston-Bennequin
vaut $\tb(L_1) = n - \w(L_0) - 1 \le -1$. La somme connexe de $L_0$ et $L_1$ est
une courbe legendrienne isotope \`a la fibre et d'enroulement~$n$.

\smallskip

\alinea{b)}
Soit $B \subset V$ une boule contenant un disque vrill\'e de~$\xi$ et soit $L_0$
une courbe legendrienne isotope \`a la fibre et disjointe de~$B$. Pour tout $n
\ge 0$, il existe dans $B$ un n\oe ud legendrien topologiquement trivial dont
l'invariant de Thurston-Bennequin vaut~$n$. La somme connexe de $L_0$ avec ce
n\oe ud fournit une courbe legendrienne isotope \`a la fibre dont l'enroulement
vaut $\w(L_0) + n + 1$. Par suite, $\w(\xi) = +\infty$.

\smallskip

\alinea{c)}
Soit $\wt\xi$ la structure induite par $\xi$ sur le rev\^etement universel $\wt
V$ de~$V$. Si $S$ n'est pas une sph\`ere, $\wt V$ est diff\'eomorphe \`a $\R^3$.
Ainsi, \`a conjugaison pr\`es, $\wt\xi$ est une  structure de contact sur $\R^2
\times \R$ transversale aux droites $\{*\} \times \R$ et invariante par les
translations verticales enti\`eres. Comme pour la proposition~\ref{p:wood}-b,
on construit sur tout domaine $[-a-1,a+1]^2 \times \R$, $a>0$, des coordonn\'ees
$(x,y,t)$ dans lesquelles la fibration est la projection $(x,y,t) \mapsto (x,y)$
et $\wt\xi$ a pour \'equation $dt - u(x,y,t)\, dx = 0$, o\`u $u(x,-a,t) = 0$ et
$u(x,y,t+1) = u(x,y,t)$ quels que soient $(x,y,t)$.

Pour tout entier~$n>0$, l'immersion $\phi_n \from [-a,a]^2 \times \R \to \R^3$,
d\'efinie --~en coordonn\'ees cylindriques au but~-- par
$$ (x,y,t) \longmapsto \Bigl( r = u(x,y,t)^{-1/2}, \;
   \theta = \frac {2\pi t} n, \; z = \frac {2\pi x} n \Bigr), $$
plonge $[-a,a]^2 \times \R/n\Z$ dans $\R^3$ priv\'e de l'axe des~$z$ et envoie
$\wt\xi$ sur la structure d'\'equation $dz = r^2\, d\theta$. Par suite, $\wt\xi$
est tendue.

Soit maintenant $L$ une courbe legendrienne isotope \`a la fibre et $\wt L$ une
pr\'eimage de~$L$ dans $\R^2 \times \R/\Z$. \'Etant donn\'e $a>0$ assez grand
pour que $[-a,a]^2 \times \R/\Z$ contienne $\wt L$, le plongement induit par
$\phi_1$ sur $[-a,a]^2 \times \R/\Z$ envoie les fibres $\{*\} \times \R/\Z$
sur des courbes deux \`a deux non enlac\'ees mais qui enlacent une fois l'axe
des~$z$. La courbe legendrienne $\phi_1(\wt L)$ est alors non nou\'ee et son
invariant de Thurston-Bennequin n'est autre que $\w(\wt L)$. Il r\'esulte donc
de l'in\'egalit\'e de Bennequin~\cite{Be} que $\w(\wt L) = \w(L) \le -1$.

Si $S$ est une sph\`ere, $\wt V$ est diff\'eomorphe \`a la sph\`ere~$\S^3$ car
$\S^2 \times \S^1$ ne porte aucune structure de contact transversale aux fibres
(th\'eor\`eme~\ref{t:milnorwood}). Ainsi, \`a conjugaison pr\`es, $\wt\xi$ est
une structure de contact sur~$\S^3$ transversale \`a la fibration de Hopf. Si on
regarde $\S^3$ comme le bord de la boule unit\'e dans~$\C^2$ -- les cercles de
Hopf \'etant les traces des droites complexes passant par~$0$ --, la forme
symplectique usuelle de~$\C^2$ est positive sur $\wt\xi$, ce qui entra\^{\i}ne
que~$\wt\xi$ est tendue~\cite{El:filling}.

Soit enfin $L$ une courbe legendrienne isotope \`a la fibre et $\wt L$ une
pr\'eimage de~$L$ dans $\wt V \simeq \S^3$. La courbe $\wt L$ est non nou\'ee
et, comme $\wt\xi$ est tendue, l'in\'egalit\'e de Bennequin assure que
$\tb(\wt L) \le -1$, donc que $\w(\wt L) = \tb(L) - 1 \le -2$. Or $\w(\wt L) =
\abs{ \Chi(V,S) } \w(L)$, donc $\w(L)$ est strictement n\'egatif.
\qed

\subsection{Le cas des fibr\'es sur la sph\`ere}

On d\'emontre ici la proposition~\ref{p:thurston} lorsque $S$ est une sph\`ere.
L'ingr\'edient cl\'e est la classification des structures de contact tendues sur
les espaces lenticulaires~\cite[th\'eo\-r\`eme~1.1]{Gi:bifurcations}. Le point
qui intervient ici est le suivant : \textit{tout espace lenticulaire porte une
seule structure de contact universellement tendue, \`a isotopie pr\`es}\footnote
{Il s'agit ici de structures de contact orientables mais pas orient\'ees.}.

\smallskip

Si $\Chi(V,S)=0$, la vari\'et\'e $V$ est diff\'eomorphe \`a $\S^2 \times \S^1$.
Or, d'apr\`es \cite{El:martinet}, $\S^2 \times \S^1$ porte une unique structure
de contact tendue qu'on peut voir par exemple comme le champ~$\xi_0$ des droites
complexes tangentes au bord de $X_\eps$, o\`u $X_\eps \subset \C^2$ est le tube
de rayon $\eps<1$ autour du cercle unit\'e $\S^1 \times \{0\}$. On observe alors
que le tore $T = \S^1 \times \eps\S^1$ est contenu dans $\partial X_\eps$ et que
$\xi_0$ est parall\`ele \`a $\C \times \{0\}$ le long de~$T$. Par suite, chaque
cercle $\S^1 \times \{w\}$, $w \in \eps\S^1$, est legendrien et d'enroulement
nul. Aucune structure de contact tendue sur~$V$ n'a donc un enroulement
strictement n\'egatif.

\smallskip

Si $\Chi(V,S) \ne 0$, le rev\^etement universel de~$V$ est diff\'eomorphe \`a
$\S^3$ et on peut voir $V$ comme suit. On regarde $\S^3$ comme le bord orient\'e
de la boule unit\'e dans $\C^2$ et, pour tout $t \in \R/\Z$, on note $\phi_t^+$
et $\phi_t^-$ les transformations de Hopf d\'efinies par
$$ \left\{ \begin{aligned}
\phi_t^+ (z,w) & = \bigl( e^{2i\pi t} z, e^{+2i\pi t} w \bigr), \\
\phi_t^- (z,w) & = \bigl( e^{2i\pi t} z, e^{-2i\pi t} w \bigr),
\end{aligned} \right. \qquad (z,w) \in \S^3 \,. $$
La fibration $\S^3 \to \S^2$ associ\'ee au flot $\phi_t^+$ (resp. $\phi_t^-$) a
pour nombre d'Euler $-1$ (resp. $+1$). Pour tout entier $n>0$, le quotient de
$\S^3$ par $\phi_{1/n}^+$ (resp. $\phi_{1/n}^-$) a donc pour nombre d'Euler~$-n$
(resp.~$n$). En outre, ce quotient n'est autre que l'espace lenticulaire
$\L_{n,1}$ (resp. $\L_{n,n-1}$) et est diff\'eomorphe \`a $V$ si $\Chi(V,S)=-n$
(resp. si $\Chi(V,S)=n$).

Soit maintenant $\xi_0$ la structure de contact usuelle sur $\S^3$, \emph{i.e.}
le champ des droites complexes tangentes \`a~$\S^3$. Comme chaque transformation
$\phi_t^\pm$ est la restriction d'une application lin\'eaire de~$\C^2$, elle
pr\'eserve $\xi_0$. Ainsi, $\xi_0$ induit sur chacun des espaces $\L_{n,1}$ et
$\L_{n,n-1}$ une structure de contact~$\xi$ qui est universellement tendue (car
$\xi_0$ est tendue d'apr\`es le th\'eor\`eme de Bennequin). D'autre part, les
orbites du flot $\phi_t^+$ sont transversales \`a~$\xi_0$ de sorte que $\xi$ est
transversale aux fibres sur $\L_{n,1}$. En revanche, sur le tore invariant
$$ \bigl\{ (z,w) \in \S^3 \mid \abs z = \abs w \bigr\}, $$
les orbites du flot $\phi_t^-$ sont tangentes \`a~$\xi_0$ et d'enroulement nul.
Par suite, $\xi$ n'a pas un enroulement strictement n\'egatif. Ces observations
terminent la d\'emonstration puisque $\L_{n,1}$ et $\L_{n,n-1}$ portent chacun
une unique structure de contact universellement tendue, \`a isotopie pr\`es
\cite[th\'eor\`eme~1.1 et lemme~4.1]{Gi:bifurcations}.
\qed

\subsection{Surfaces convexes}

On introduit ici quelques notions et r\'esultats techniques qui seront utiles
dans les d\'emonstrations \`a venir.

\begin{definition} \label{d:ctc}
Soit $F$ une surface, orientable et compacte, plong\'ee dans une vari\'et\'e de
contact $(M,\xi)$ de dimension~$3$. On dit que $F$ est \emph{convexe} si elle
admet un voisinage tubulaire trivialis\'e $U \simeq F \times \R$ dans lequel les
translations verticales pr\'eservent~$\xi$. Un tel voisinage~$U$ sera dit \emph
{homog\`ene}.

La convexit\'e de~$F$ ne d\'epend que du germe de~$\xi$ le long de~$F$, donc du
feuilletage caract\'eristique~$\xi F$ de~$F$ form\'e des courbes int\'egrales du
champ de droites $\xi \cap TF$. Si $F$ est close, elle se traduit explicitement
comme suit. On dit qu'une multi-courbe $\Gamma \subset F$ --~union finie de
courbes ferm\'ees, simples et disjointes~--
\emph{scinde} $\xi F$ si, sur la surface compacte \`a bord~$F_\Gamma$ obtenue en
d\'ecoupant~$F$ le long de~$\Gamma$, le feuilletage induit par~$\xi F$ est
port\'e par un champ de vecteurs qui sort sur $\partial F_\Gamma$ et qui dilate
l'aire (\emph{i.e.} une certaine forme d'aire sur~$F_\Gamma$). Lorsqu'une telle
multi-courbe existe, elle est unique \`a isotopie pr\`es parmi les multi-courbes
qui scindent.
\end{definition}

Les r\'esultats de \cite{Gi:convexite} montrent qu'une surface close $F \subset
(M,\xi)$ est convexe si et seulement si son feuilletage $\xi F$ est scind\'e. En
particulier, si $U$ est un voisinage tubulaire et homog\`ene d'une surface
convexe~$F$, l'ensemble not\'e~$\Gamma_{\!U}$ des points de~$F$ o\`u $\xi$ est
tangente aux fibres de~$U$ est une multi-courbe qui scinde $\xi F$.

Cette caract\'erisation permet de montrer que les surfaces closes convexes sont
g\'en\'e\-riques. En outre, elles sont tr\`es maniables et certaines ont un
r\'eel int\'er\^et g\'eom\'etrique :

\begin{exemple}
Suivant une suggestion de V.~I.~{\sc Arnold}, on appellera \emph{surface
clairaldienne} toute surface compacte convexe $F \subset (M,\xi)$ qui est munie
d'une fibration en cercles legendriens, au-dessus de l'intervalle ou du cercle.
Topologiquement, une telle surface est donc un anneau, un tore ou une bouteille
de Klein.

Si $\pi \from V_0 \to S_0$ est une fibration legendrienne, l'image inverse
$\pi^{-1}(C)$ de toute courbe simple $C \subset S_0$ (ferm\'ee ou non) est une
surface clairaldienne. En effet, tout flot local transversal \`a~$C$ dans $S_0$
se rel\`eve naturellement dans~$V_0$ en un flot de contact transversal \`a
$\pi^{-1}(C)$. Inversement, toute surface clairaldienne~$F \subset (M,\xi)$ est
localement de ce type $\pi^{-1}(C)$, o\`u $V_0$ est un voisinage homog\`ene
quelconque de~$F$ et $S_0$ un fibr\'e en intervalles au-dessus de~$C$. Ainsi,
toute surface clairaldienne orientable --~diff\'eomorphe \`a $\S^1 \times C$
o\`u $C$ est l'intervalle ou le cercle~-- poss\`ede un voisinage tubulaire
homog\`ene $U \simeq \S^1 \times C \times \R$ dans lequel $F = \S^1 \times C
\times \{0\}$ et $\xi$ a une \'equation de la forme
$$ \cos (n\pi x) \, dy - \sin (n\pi x) \, dt = 0, \qquad
   (x,y,t) \in \S^1 \times C \times \R, $$
o\`u $n$ est un entier strictement positif, pair d\`es que $\xi$ est orientable.
Le nombre $-n/2$ n'est autre que l'enroulement de la structure autour des
cercles legendriens qui fibrent~$F$.
\end{exemple}

\smallskip

Le lemme qui suit est une version relative d'un r\'esultat \'etabli dans~\cite
{Gi:convexite}.

\begin{lemme} \label{l:ctc}
Soit $F \subset (M,\xi)$ une surface convexe, $U$ un voisinage tubulaire de~$F$
homog\`ene, $P$ un compact de $M$ dont l'intersection avec~$F$ est satur\'ee par
$\xi F$ et $\sigma$ un feuilletage de~$F$ scind\'e par $\Gamma_{\!U}$ et \'egal
\`a $\xi F$ pr\`es de $P \cap F$. Il existe alors une isotopie de plongements
$\phi_t \from F \to U$, $t \in [0,1]$, ayant les propri\'et\'es suivantes\up{ :}
\begin{enumerate}
\item[\up{1)}]
$\phi_0$ est l'inclusion\up{ ;}
\item[\up{2)}]
pour tout $t \in [0,1]$, la surface $\phi_t(F)$ est transversale aux fibres
de~$U$\up{ ;}
\item[\up{3)}]
le feuilletage caract\'eristique $\xi\, \phi_1(F)$ n'est autre que $(\phi_1)_*
\sigma$\up{ ;}
\item[\up{4)}]
l'intersection $\phi_t(F) \cap P$ co\"{\i}ncide avec $F \cap P$ pour tout $t \in
[0,1]$.
\end{enumerate}
\end{lemme}

\begin{proof}
La proposition II.3.6 de~\cite{Gi:convexite} donne une isotopie $\ol\phi_t$ qui
satisfait aux conditions 1--3) et laisse fixe un voisinage de $F \cap P$. Pour
obtenir~4), on remarque qu'on garde 1--3) si on compose $\ol\phi_t$ au but par
une isotopie de contact partant de l'identit\'e et pr\'eservant la structure
fibr\'ee du tube~$U$. On note alors $\delta_s \from U \simeq F \times \R \to U$,
$s \in ]0,1]$, l'homoth\'etie de rapport~$s$ dans les fibres et on observe que,
pour $s_0$ assez petit, toutes les surfaces $\delta_{s_0} \circ \ol\phi_t (F)$,
$t \in [0,1]$, coupent~$P$ exactement suivant $F \cap P$. D'autre part, comme
$\xi$ est $\R$-invariante dans~$U$, elle y admet une \'equation de Pfaff du type
$\beta + u\,dt = 0$, o\`u $t$ d\'ecrit~$\R$ et $\beta$, $u$ sont respectivement
une $1$-forme et une fonction sur~$F$. Chaque structure $\xi_s = (\delta_s)_*
\xi$ a ainsi pour \'equation $\beta + (u/s)\,dt = 0$ et est donc encore $\R
$-invariante. La m\'ethode du chemin fournit alors une isotopie $\psi_s$ de~$U$,
$s \in [s_0,1]$, qui est $\R$-\'equivariante, envoie $\xi_s$ sur~$\xi=\xi_1$ et
d\'eplace les points \emph{horizontalement} en respectant tous les feuilletages
$\xi_s (F \times \{*\}) = \xi (F \times \{*\})$. En particulier, comme
l'intersection $F \cap P$ est satur\'ee par $\xi F$, elle est pr\'eserv\'ee par
l'isotopie~$\psi_s$. On prend alors une fonction lisse $s \from [0,1] \to ]0,1]$
qui vaut~$1$ en ~$0$ mais devient vite tr\`es petite. L'isotopie
$$ \phi_t = \psi_{s(t)} \circ \delta_{s(t)} \circ \ol\phi_t
   \from F \longrightarrow U $$
v\'erifie toutes les propri\'et\'es voulues.
\end{proof}

\subsection{Redressement des tores}

Un ingr\'edient cl\'e dans la th\`ese de W.~{\sc Thurston}~\cite{Th:these} est
le r\'esultat suivant, d\^u ind\'ependamment \`a R.~{\sc Roussarie}~\cite{Ro} :
Dans une vari\'et\'e de dimension~$3$ munie d'un feuilletage de codimension~$1$
sans composantes de Reeb, tout tore incompressible plong\'e est isotope \`a une
feuille ou \`a un tore transversal au feuilletage. Dans les vari\'et\'es de
contact, les techniques de~\cite{Gi:convexite} permettent d'\'etablir un fait
analogue tr\`es utile pour d\'emontrer le th\'eor\`eme~\ref{t:thurston} :

\begin{lemme} \label{l:roussarie}
Soit $\xi$ une structure de contact sur~$V$ et $R \subset S$ une sous-surface
compacte, connexe et \`a bord non vide. Si l'enroulement $\w(\xi)$ de~$\xi$ est
strictement n\'egatif, $\xi$ est isotope \`a une structure de contact~$\xi'$
pour laquelle, au-dessus de~$R$, toutes les fibres sont legendriennes et ont un
enroulement \'egal \`a $\w(\xi)$.
\end{lemme}

\begin{proof}
On regarde la surface~$R$ comme un voisinage r\'egulier d'un bouquet de cercles
$K = \bigvee_{i=1}^k K_i$ dans $S$ et on note~$q$ le sommet de~$K$. Quitte \`a
faire une premi\`ere isotopie, on suppose que la fibre~$L$ au-dessus de~$q$ est
legendrienne et que son enroulement vaut $\w(\xi)$. Pour tout entier~$n>0$, on
peut trouver, sur un voisinage tubulaire $W$ de~$L$, des coordonn\'ees $(x,y,t)
\in \D^2 \times \S^1$ dans lesquelles $L = \{0\} \times \S^1$ et $\xi$ a pour
\'equation
$$ \cos (2n\pi t) \, dx - \sin (2n\pi t) \, dy = 0 \,.$$
La projection $\pi_W \from W \to \D^2$, \ $(x,y,t) \mapsto (x,y)$, est alors une
fibration legendrienne et l'enroulement de~$\xi$ autour des fibres de~$\pi_W$
vaut~$-n$. Ainsi, pour $n = -\w(\xi)$ (qui est un entier strictement positif),
$\pi_W$ induit la m\^eme trivialisation normale de~$L$ que la fibration $\pi
\from V \to S$. Il existe donc une isotopie $\phi_t \from V \to V$, $t \in
[0,1]$, qui a les propri\'et\'es suivantes :
\begin{itemize}
\item
$\phi_0 = \id$ ;
\item
$\phi_t(L) = L$ pour tout $t \in [0,1]$ ;
\item
$\phi_1$ envoie chaque fibre de $\pi_W$ sur une fibre de~$\pi$.
\end{itemize}
Quitte \`a remplacer~$\xi$ par $(\phi_1)_*\xi$, on suppose d\'esormais que~$\xi$
est tangente aux fibres de~$V$ au-dessus d'un voisinage compact~$Q$ de~$q$. On
r\'eduit~$Q$ au besoin pour que $Q \cap K$ soit connexe. Une version relative
facile des r\'esultats de \cite{Gi:convexite} permet alors de d\'eformer~$\xi$
par une $\classe\infty$-petite isotopie relative \`a $M = \pi^{-1}(Q)$ pour que
chaque anneau $\pi^{-1}(K_i \setminus Q)$, $1 \le i \le k$, soit convexe. On
lisse ensuite chaque~$K_i$ dans~$Q$ en une courbe~$K'_i$. Les tores $F_i =
\pi^{-1}(K'_i)$ sont alors convexes. On en prend des voisinages homog\`enes
respectifs~$U_i$ et on pose $\Gamma_{\!i} = \Gamma_{\!U_i}$ (cf. d\'efinition
\ref{d:ctc}).

\begin{assertion}
L'intersection g\'eom\'etrique $\Card (L \cap \Gamma_{\!i})$ de~$L$ avec~$
\Gamma_{\!i}$ est \'egale au module $\abs{ [L] \cdot [\Gamma_{\!i}] }$ de leur
intersection alg\'ebrique --~toutes les composantes connexes de~$\Gamma_{\!i}$
\'etant orient\'ees dans le m\^eme sens.
\end{assertion}

\begin{proof}[Preuve]
L'enroulement $\w(L)$ se lit sur~$F_i$ comme $-(1/2) \Card (L\cap\Gamma_{\!i})$.
Si l'assertion est fausse, $F_i$ porte une courbe ferm\'ee simple~$C$ isotope
\`a~$L$ qui intersecte~$\Gamma_{\!i}$ moins que~$L$ (g\'eom\'etriquement). On
peut alors construire sans peine sur~$F_i$ un feuilletage singulier~$\sigma$ qui
est scind\'e par~$\Gamma_{\!i}$ et pour lequel la courbe~$C$ est satur\'ee (voir
\cite[exemple II.3.7]{Gi:convexite}). Le lemme~II.3.6 de~\cite{Gi:convexite}
(version absolue du lemme~\ref{l:ctc}) fournit alors un plongement $\phi$ de
$F_i$ dans $U_i$ --~isotope \`a l'inclusion~-- dont l'image a pour feuilletage
caract\'eristique $\phi_*\sigma$. L'enroulement de~$\xi$ autour de la courbe
legendrienne $L' = \phi(C)$ vaut alors
$$ \w(L') = -\tfrac12 \Card (C \cap \Gamma_{\!i})
 > -\tfrac12 \Card (L \cap \Gamma_i) = \w(L), $$
ce qui contredit le fait que $\w(L)$ est \'egal \`a $\w(\xi)$.
\end{proof}

\smallskip

L'assertion ci-dessus permet de d\'eformer $\xi$, par une isotopie relative \`a
$M$ laissant les tores~$F_i$ invariants, de telle sorte que chaque fibre de $\pi
\res {F_i}$ ait, avec~$\Gamma_{\!i}$, une intersection g\'eom\'etrique \'egale
au module de son intersection alg\'ebrique. Cette condition \'etant remplie, il
existe sur~$F_i$ un feuilletage singulier~$\sigma_i$ ayant les propri\'et\'es
suivantes :
\begin{itemize}
\item
$\sigma_i$ est scind\'e par~$\Gamma_{\!i}$ ;
\item
$\sigma_i$ co\"{\i}ncide avec $\xi F_i$ dans $F_i \cap M$ ;
\item
chaque fibre de $\pi \res {F_i}$ est satur\'ee par~$\sigma_i$.
\end{itemize}
Le lemme~\ref{l:ctc} donne alors, pour $1 \le i\le k$, un plongement $\phi_i$ de
$F_i$ dans~$U_i$ --~isotope \`a l'inclusion~-- dont l'image a pour feuilletage
caract\'eristique $(\phi_i)_*\sigma_i$ et a m\^eme intersection que $F_i$ avec
le compact
$$ P_i = M \cup \phi_1(F_1) \cup \dots \cup \phi_{i-1}(F_{i-1}) \cup
   F_{i+1} \cup \dots \cup F_k \,. $$
Il existe donc un diff\'eomorphisme~$\phi$ de~$V$, isotope \`a l'identit\'e, qui
\guil{prolonge} simultan\'ement tous les plongements~$\phi_i$. Par construction,
la structure de contact~$\phi^*\xi$ imprime le feuilletage $\sigma_i$ sur chaque
tore~$F_i$ et est ainsi tangente aux fibres de~$\pi$ au-dessus de $Q \cup K$.
Comme tous les tores~$F_i$ sont convexes (voire clairaldiens), il est facile de
rendre la fibration~$\pi$ legendrienne au-dessus de~$R$ par une ultime isotopie
relative \`a $\pi^{-1}(Q \cup K)$.
\end{proof}

\subsection{Structures de contact sur le tore plein}

La d\'emonstration de la proposition~\ref{p:thurston} passe par une analyse des
structures de contact tendues sur le tore plein $W = \D^2 \times \S^1$. Cette
analyse est men\'ee dans~\cite{Gi:bifurcations} et on en pr\'esente ici quelques
conclusions utiles. Pour cela, on rappelle qu'un feuilletage du tore $\T^2$ est
une \emph{suspension} s'il est non singulier et si toutes ses feuilles coupent
une m\^eme courbe transversale ferm\'ee, simple et connexe. D'autre part, on
observe qu'une structure de contact sur $\D^2 \times \S^1$ est orientable si et
seulement si elle l'est le long des fibres de la projection $\D^2 \times \S^1
\to \D^2$.

\begin{lemme} \label{l:bfs}
Soit $\xi$ une structure de contact orientable et tendue sur $W = \D^2 \times
\S^1$. On suppose que le feuilletage caract\'eristique $\xi \, \partial W$ est
scind\'e par une multi-courbe ayant $2n$~composantes connexes et que ses
\'eventuelles
singularit\'es forment des cercles lisses. Il existe alors $n$~anneaux disjoints
$A_i$ plong\'es dans $W$ et ayant les propri\'et\'es suivantes\up{ :}
\begin{itemize}
\item
chaque composante de $\partial A_i$ est une courbe de singularit\'es ou une
feuille ferm\'ee de $\xi\, \partial W$\up{ ;}
\item
chaque feuilletage $\xi A_i$ est constitu\'e de cercles parall\`eles au bord.
\end{itemize}
\end{lemme}

\begin{proof}
Si le feuilletage $\xi\,\partial W$ est une suspension, les anneaux~$A_i$ sont
directement fournis par la proposition 3.15 de~\cite{Gi:bifurcations} : dans la
terminologie de cet article, ce sont les anneaux du feuillage d'une structure de
contact \'el\'ementaire isotope \`a~$\xi$ relativement au bord. On va maintenant
adapter l'argument au cas o\`u $\xi\, \partial W$ est un feuilletage scind\'e
dont les singularit\'es forment des cercles. Dans ce cas, l'\'etude des surfaces
convexes (voir les sections 2.B et 2.C de \cite{Gi:bifurcations}) montre que
$\xi$ est isotope, relativement au bord, \`a une structure de contact $\xi'$
ayant les propri\'et\'es suivantes :
\begin{itemize}
\item
chaque tore $T_a = a\S^1 \times \S^1$, $1/2 \le a \le 1$, est convexe dans
$(W,\xi')$ ;
\item
le feuilletage $\xi T_{1/2}$ est une suspension ;
\item
les singularit\'es \'eventuelles de chaque feuilletage $\xi' T_a$, $1/2 < a \le
1$, forment des cercles.
\end{itemize}
Dans $W' = (1/2)\D^2 \times \S^1$, la proposition~3.15 de \cite{Gi:bifurcations}
donne, comme avant, des anneaux~$A'_i$ qui conviennent pour la restriction de
$\xi'$. D'autre part, dans $W \setminus \Int W'$, l'union des feuilles ferm\'ees
et des singularit\'es de tous les feuilletages $\xi' T_a$, $1/2 \le a \le 1$,
forme $2n$~anneaux disjoints qui compl\`etent les $A'_i$ en les anneaux $A_i$
cherch\'es.
\end{proof}

\begin{proposition} \label{p:bfs}
Soit $\xi$ une structure de contact orientable et universellement tendue sur $W
= \D^2 \times \S^1$. On suppose que le feuilletage caract\'eristique $\xi \,
\partial W$ est scind\'e par deux courbes et que ses singularit\'es forment deux
cercles lisses. La structure $\xi$ est alors isotope, relativement au bord, \`a
une structure de contact qui est transversale \`a $\{0\} \times \S^1$ et imprime
une suspension sur chaque tore $a\S^1 \times \S^1$, $0 < a < 1$.
\end{proposition}

\begin{proof}
Soit $\sigma$ un feuilletage de $\partial W$ scind\'e par deux courbes et dont
les singularit\'es forment deux cercles lisses.
D'apr\`es le
th\'eor\`eme 1.6 de \cite{Gi:bifurcations}, compl\'et\'e par le lemme 3.13, les
structures de contact universellement tendues sur~$W$ qui impriment~$\sigma$ sur
$\partial W$ forment au plus deux classes d'isotopie relative au bord. Chaque
classe est caract\'eris\'ee par la classe d'isotopie de l'anneau que fournit le
lemme~\ref{l:bfs}. Autrement dit, il y a deux classes (resp. une) s'il y a dans
$W$, \`a isotopie relative au bord pr\`es, deux anneaux (resp. un seul) qui
s'appuient sur les cercles singuliers de $\sigma$. On exhibe ci-dessous des
structures de contact universellement tendues explicites dans chaque classe et
on constate qu'elles satisfont les propri\'et\'es requises.

Sur $\R^2 \times \S^1$ muni de coordonn\'ees cylindriques $(r,\theta,z)$, $z \in
\R/2\pi\Z$, l'\'equation de Pfaff $(1-r^4)\, dz + r^2\, d\theta = 0$ d\'efinit
une structure de contact~$\zeta$ universellement tendue (l'\'equation de~$\zeta$
d\'efinit sur $\R^3$ la structure de contact ordinaire). De plus, pour tout $r>0
$, le feuilletage caract\'eristique du tore de rayon~$r$ autour de $\{0\} \times
\S^1$ est le feuilletage lin\'eaire de pente $dz/d\theta = r^2 \!/ (r^4-1)$.
\'Etant donn\'e des entiers $p$ et $q$ premiers entre eux, $q>0$, il existe donc
un unique r\'eel $r = r(p,q)$ tel que les caract\'eristiques du tore $\partial
(r\D^2) \times \S^1$ aient pour classe d'homologie $(p,q)$. On consid\`ere alors
les plongements $\psi_{p,q}^\pm \from W' = (1/2)\D^2 \times \S^1 \to \R^2 \times
\S^1$ donn\'es par
$$ \psi_{p,q}^\pm (ae^{is}, t) =
   \Bigl( 2ar(p,q) \bigl( 1 \pm \tfrac a q \cos(qs-pt) \bigr), \,
   s + \tfrac a q \sin \bigl( 2(qs-pt) \bigr), \, t \Bigr) \,. $$
Les structures de contact induites, $\zeta_\pm = (\psi_{p,q}^\pm)^*\zeta$, sont
universellement tendues, transversales \`a $\{0\} \times \S^1$  et impriment une
suspension sur chaque tore $T_a = a\S^1 \times \S^1$, $0 < a < 1/2$. En outre,
les feuilletages $\sigma_\pm = \zeta_\pm\, \partial W'$ sont tous deux scind\'es
par deux courbes et ont deux cercles de singularit\'es qui sont communs et qu'on
note $C_0$, $C_1$. On observe d'autre part que l'image inverse du tore de rayon
$r(p,q)$ par $\psi_{p,q}^\pm$ est un anneau $A_\pm$ qui a les propri\'et\'es
d\'ecrites au lemme~\ref{l:bfs}. De plus, $A_-$ n'est isotope \`a $A_+$
relativement \`a son bord que si $q=1$.

\begin{figure}[ht]
\centering
\includegraphics[width=.8\textwidth]{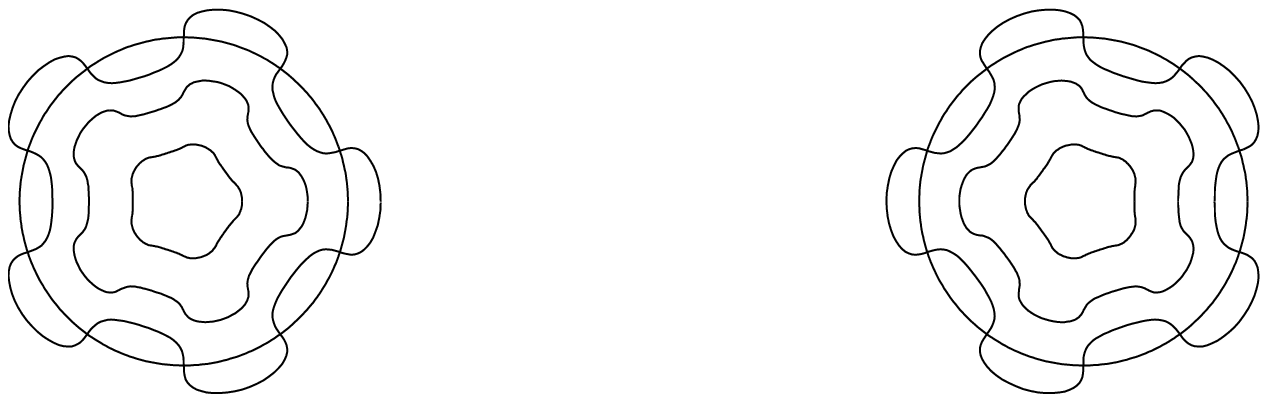}
\legende{Pour $q=5$, les images respectives par $\psi_{p,q}^+$ et $\psi_{p,q}^-$
des cercles de rayons $1/6$, $1/3$ et $1/2$ dans $\D^2 \times \{0\}$.}
\end{figure}

On suppose maintenant que $(p,q)$ est la classe des cercles singuliers de~$
\sigma$ et on choisit dans $W \setminus \Int W'$ deux anneaux $B_0$, $B_1$ qui
sont transversaux aux tores $T_a = a\S^1 \times \S^1$, $1/2 \le a \le 1$, et qui
s'appuient d'un c\^ot\'e sur les cercles singuliers de $\sigma$, de l'autre sur
$C_0$ et $C_1$. Il existe alors sur $W \setminus \Int W'$  deux structures de
contact $\eta_-$ et $\eta_+$ satisfaisant aux conditions suivantes (voir
\cite[lemme 2.3]{Gi:bifurcations}) :
\begin{itemize}
\item
$\eta_\pm\, \partial W' = \sigma_\pm$ et $\eta_\pm \partial W = \sigma$ ;
\item
chaque feuilletage $\eta_\pm T_a$, $1/2 \le a \le 1$, est scind\'e et a deux
cercles de singularit\'es, \`a savoir les cercles $B_i \cap T_a$, $i=0,1$.
\end{itemize}
Les structures de contact $\xi_\pm = \zeta_\pm \cup^{}_{\partial W'} \eta_\pm$
sont universellement tendues et ne sont isotopes, relativement au bord, que si
les anneaux $A_\pm \cup B_0 \cup B_1$ le sont (\emph{i.e.} si $q=1$). En outre,
elles sont transversales \`a $\{0\} \times \S^1$, impriment $\sigma$ sur le bord
$\partial W$ et une suspension sur $T_a$ pour tout $a \in \op]0,1/2\cl[$. Sans
d\'etruire ces propri\'et\'es, une $\classe\infty$-petite isotopie convenable
--~\`a support dans un voisinage de $(B_0 \cup B_1) \cap \Int W$~-- permet de
perturber $\xi_\pm$ en une structure de contact qui imprime une suspension sur
tous les tores $T_a$, $0 < a < 1$.
\end{proof}

\subsection{Mise en position transversale}

On d\'emontre ici la proposition~\ref{p:thurston} lorsque $S$ n'est pas une
sph\`ere. Dans un premier temps, $\xi$ d\'esigne juste une structure de contact
d'enroulement strictement n\'egatif sur~$V$. On note~$g\ge1$ le genre de~$S$ et
$K$ un bouquet de $2g$~cercles sur~$S$ ayant pour compl\'ementaire un disque.
Compte tenu du lemme~\ref{l:roussarie}, on suppose que, au-dessus d'un voisinage
compact r\'egulier~$R$ de~$K$, les fibres de~$\pi$ sont legendriennes et ont
pour enroulement $\w(\xi)$. On note $D$ le disque ferm\'e $S \setminus \Int R$
et $W$ le tore plein $\pi^{-1}(D)$ qu'on param\`etre par $\D^2 \times \S^1$ de
telle sorte que la fibration $\pi \res W$ soit la projection sur le premier
facteur. Par construction, $\partial W$ est un tore clairaldien.

\begin{lemme} \label{l:deux}
La multi-courbe~$\Gamma$ qui scinde le feuilletage $\xi\, \partial W$ a deux
composantes connexes.
\end{lemme}

\begin{proof}
Soit $2n$ le nombre (pair) de composantes connexes de~$\Gamma$. D'apr\`es le
lemme~\ref{l:bfs}, les $2n$~courbes de singularit\'es de $\xi\, \partial W$
bordent $n$~anneaux~$A_i$ disjoints et plong\'es dans~$W$ dont les feuilletages
caract\'eristiques $\xi A_i$ sont form\'es de cercles parall\`eles au bord. On
indexe les~$A_i$ de telle sorte que $A_1$ soit ext\'erieurissime, c'est-\`a-dire
d\'ecoupe~$W$ en deux tores pleins dont l'un, not\'e~$W_1$, se r\'etracte par
d\'eformation sur~$A_1$ et ne contient aucun $A_i$, $i>1$. On prend ensuite un
voisinage collier $N_1 \cong A_1 \times [0,1]$ de $A_1 = A_1 \times \{0\}$ dans
$\Adh (W \setminus W_1)$ dont le bord lat\'eral $\partial A_1 \times [0,1]$ est
inclus dans $\partial W$. Pour tout $s\ne0$ assez petit, les caract\'eristiques
de l'anneau $A_1 \times \{s\}$ vont d'un bord \`a l'autre. On d\'esigne alors
par~$W'$ un tore plein obtenu en arrondissant les angles de
$$ W \setminus \Bigl( W_1 \cup \bigl( A_1 \times [0,s \mathclose[ \bigr) \Bigr),
   \qquad \text{$s>0$ petit.} $$
Ainsi construit, $W'$ est un tore plein isotope \`a~$W$ et on note $\psi \from
\partial W \to V$ un plongement isotope \`a l'inclusion et dont l'image est le
tore $\partial W'$.

Si $n>1$, le feuilletage $\xi\, \partial W'$ a $2n-2$ courbes de singularit\'es
et aucune feuille r\'eguli\`ere ferm\'ee. Par suite, $\partial W'$ est convexe
et la multi-courbe $\Gamma'$ qui scinde $\xi\, \partial W'$ compte
$2n-2$~composantes
connexes, toutes isotopes aux composantes de $\psi(\Gamma)$. En notant $L$ une
fibre de $\pi$ dans $\partial W$, il existe sur $\partial W'$ une courbe~$L'$
isotope \`a $\psi(L)$ et v\'erifiant
$$ \Card (L' \cap \Gamma') < 2 \abs{ \w(\xi) } = \Card (L \cap \Gamma) \,. $$
Or le lemme~\ref{l:ctc} fournit une isotopie de plongements $\phi_s \from
\partial W' \to V$, $s \in [0,1]$, ayant les propri\'et\'es suivantes :
\begin{itemize}
\item
chaque feuilletage $\xi\, \phi_s(\partial W')$, $s \in [0,1]$, est scind\'e par
$\phi_s(\Gamma')$ ;
\item
la courbe $\phi_1(L')$ est legendrienne.
\end{itemize}
Comme dans la preuve du lemme~\ref{l:roussarie}, l'enroulement de $\xi$ autour
de $\phi_1(L')$ vaut alors 
$$ \w \bigl( \phi_1(L') \bigr) = - \tfrac 12 \Card (L' \cap \Gamma')
 > \w(\xi), $$
ce qui contredit la d\'efinition de $\w(\xi)$.
\end{proof}

\smallskip

Le second ingr\'edient dans la preuve de la proposition~\ref{p:thurston} est le
lemme suivant :

\begin{lemme} \label{l:univ}
Si la structure de contact~$\xi$ est universellement tendue, sa restriction au
tore plein $W = \pi^{-1}(D)$ l'est aussi.
\end{lemme}

\begin{proof}
Soit $n\ge1$ un entier quelconque. Comme la surface~$S$ n'est pas une sph\`ere,
elle poss\`ede un rev\^etement connexe \`a $n$~feuillets~$S_n$. Le rappel $V_n
\to S_n$ du fibr\'e $V \to S$ au-dessus de~$S_n$ a alors pour nombre d'Euler
$\Chi(V_n,S_n) = n \Chi(V,S)$. Par suite, il existe un rev\^etement fibr\'e \`a
$n$~feuillets $\wt V_n \to V_n$ qui, au-dessus de chaque fibre de la projection
$V_n \to S_n$, induit un rev\^etement cyclique non trivial du cercle. L'image
inverse de~$W$ dans $\wt V_n$ est la r\'eunion disjointe de $n$~tores pleins et
chacun d'eux se projette sur~$W$ par un rev\^etement de degr\'e~$n$. Comme $\xi$
se rel\`eve sur $\wt V_n$ en une structure de contact tendue, sa restriction \`a
$W$ induit une structure de contact tendue sur tout rev\^etement fini de~$W$,
donc aussi sur le rev\^etement universel.
\end{proof}

\smallskip

On termine maintenant la d\'emonstration de la proposition~\ref{p:thurston}.
On note donc $\xi$ une structure de contact sur $V$ qui est orientable le long
des fibres, universellement tendue et d'enroulement strictement n\'egatif. En
outre, compte tenu des lemmes \ref{l:roussarie}, \ref{l:deux}, \ref{l:univ} et
de la proposition~\ref{p:bfs}, on suppose que $\xi$ satisfait les propri\'et\'es
suivantes :
\begin{itemize}
\item
toutes les fibres au-dessus du voisinage $R$ de $K$ sont legendriennes ;
\item
dans~$W = \pi^{-1}(D) \cong \D^2 \times \S^1$, la structure de contact~$\xi$ est
transversale \`a $\{0\} \times \S^1$ et imprime une suspension sur chaque tore
$T_a = a\S^1 \times \S^1$, $0 < a < 1$.
\end{itemize}
Comme l'enroulement $\w(\xi)$ est strictement n\'egatif, les feuilletages~$\xi
T_a$, $0 < a < 1$, n'ont aucune feuille ferm\'ee isotope \`a la fibre et sont
par suite tous isotopes \`a des feuilletages transversaux aux fibres. Quitte \`a
d\'eformer~$\xi$ par une isotopie \`a support dans~$W$, on peut donc supposer
que toutes les fibres, au-dessus de $\Int D$, sont transversales \`a~$\xi$. On
oriente alors $\xi$ sur un voisinage de~$W$ pour que les dites fibres $\{*\}
\times \S^1 \subset \Int W$ soient des transversales positives. On se donne par
ailleurs un champ de vecteurs legendrien~$\tau$ sur~$V$ qui est nul sur~$W$ et
transversal aux fibres sur $V \setminus W$. Si on pousse~$\xi$ par le flot de
$\tau$ pendant un bref instant, on obtient une structure de contact~$\xi'$ qui
est partout transversale aux fibres sauf le long de $\partial W$ o\`u elle reste
tangente aux fibres. De plus, si on pousse dans la bonne direction, les fibres
de part et d'autre de $\partial W$ sont transversales dans le m\^eme sens. Une
isotopie $\classe\infty$-petite permet alors de rendre $\xi'$ transversale \`a
toutes les fibres.
\qed

\subsection{Un exemple}

Pour clore cette partie, on montre que l'alternative offerte par le th\'eor\`eme
\ref{t:thurston} est optimale au sens o\`u, lorsque~$\xi$ n'est pas isotope \`a
une structure transversale aux fibres, il est parfois indispensable de passer
\`a un rev\^etement fini de~$V$ pour trouver une courbe legendrienne isotope \`a
la fibre et d'enroulement nul.

\begin{proposition} \label{p:noeuds}
Si $S$ est un tore et si $\Chi(V,S)$ est un nombre n\'egatif assez grand, $V$
porte une structure de contact virtuellement vrill\'ee d'enroulement strictement
n\'egatif.
\end{proposition}

\begin{proof}
Soit $\zeta$ la structure de contact d'\'equation $dz - y\, dx = 0$ sur $\R^3$.
On se donne un n\oe ud~$K$ transversal \`a~$\zeta$ et un voisinage tubulaire~$W$
de~$K$ dans lequel $\zeta$ a pour \'equation $dt + r^2\, d\theta= 0$, o\`u~$t$
param\`etre~$K$ et $(r,\theta)$, $r \le \eps$, sont des coordonn\'ees polaires
normales. On note $\l(K)$ l'\emph{auto-enlacement} de~$K$ dans $(\R^3,\zeta)$,
enlacement de~$K$ avec $K + \eps\partial_y$, et $K'$ une stabilisation de~$K$
dans~$W$, c'est-\`a-dire un n\oe ud topologiquement isotope \`a~$K$ dans~$W$,
transversal \`a~$\zeta$ et d'auto-enlacement $\l(K') = \l(K) - 2$ \ (voir \cite
{Be}). Le th\'eor\`eme de Darboux assure alors qu'il existe un isomorphisme~$
\phi$ de $(W,\zeta)$ sur un voisinage tubulaire $(W',\zeta)$ de~$K'$ qu'on peut
prendre aussi petit qu'on veut. Comme $\phi$ envoie les m\'eridiens de~$W$ sur
des m\'eridiens de~$W'$, la vari\'et\'e~$V$ obtenue \`a partir de $W \setminus
\Int W'$ en identifiant par~$\phi$ les deux composantes du bord est fibr\'ee en
cercles au-dessus du tore. On note~$\xi$ la structure de contact induite par
$\zeta$ sur~$V$ et on observe que les propri\'et\'es suivantes sont
satisfaites :
\begin{itemize}
\item
le nombre d'Euler de ~$V$ est n\'egatif, d'autant plus grand que $W'$ est plus
petit ;
\item
au signe pr\`es, la classe d'Euler de~$\xi$ est en dualit\'e de Poincar\'e avec
la fibre de~$V$.
\end{itemize}
Cette seconde propri\'et\'e montre imm\'ediatement que~$\xi$ n'est pas isotope
\`a une structure transversale aux fibres. Par ailleurs, si $(V,\xi)$ contenait
une courbe legendrienne isotope \`a la fibre et d'enroulement nul, celle-ci se
rel\`everait --~peut-\^etre pas dans $W \setminus W'$ mais dans $W \setminus
\phi^n(W')$ pour $n$~assez grand~-- en une courbe legendrienne de $(W,\zeta)$
bordant un disque m\'eridien vrill\'e. Or la structure $\zeta$ est tendue
d'apr\`es le th\'eor\`eme de Bennequin.
\end{proof}

\smallskip

Au prix de quelques efforts suppl\'ementaires, l'exemple ci-dessus r\'ev\`ele
aussi qu'on ne peut pas se contenter de consid\'erer des rev\^etements finis du
type $\wt V = \rho^* V$ o\`u $\rho$ est un rev\^etement fini de~$S$ : il faut
en g\'en\'eral d\'eplier les fibres.

\section{D\'enombrement des structures de contact trans\-versales}

\subsection{Comment contacter Matsumoto-Ghys}

$V$ d\'esigne toujours une vari\'et\'e connexe orient\'ee fibr\'ee en cercles
au-dessus d'une surface close~$S$. On classifie ici les structures de contact
directes et transversales aux fibres sur~$V$. Pour ce qui est des feuilletages,
aucune classification topologique g\'en\'erale n'est connue ni m\^eme attendue
(voir~\cite{Gh:euler}). Toutefois, lorsque $\Chi(V,S)$ vaut $\pm \Chi(S)$, les
travaux de S.~{\sc Matsumoto} et \'E.~{\sc Ghys} montrent que les feuilletages
$\classe2$ transversaux aux fibres sur~$V$ sont tous topologiquement conjugu\'es
\cite{Ma} et forment, \`a conjugaison diff\'erentiable pr\`es, une vari\'et\'e
hom\'eomorphe \`a l'espace de Teichm\"uller de~$S$ \cite{Gh:fuchs}. Pour les
structures de contact, le th\'eor\`eme qui suit donne, sans autre restriction
sur $\Chi(V,S)$ que l'in\'egalit\'e du th\'eor\`eme~\ref{t:milnorwood}, une
classification compl\`ete. On rappelle que, d'apr\`es les propositions~\ref
{p:calculs} et~\ref{p:thurston}, une structure de contact sur~$V$ est isotope
\`a une structure transversale aux fibres si et seulement si elle est
universellement tendue et d'enroulement strictement n\'egatif entier.

\begin{theoreme} \label{t:ghys}
Soit $V$ une vari\'et\'e connexe et orient\'ee, fibr\'ee en cercles au-dessus
d'une surface close et orientable~$S$.  On suppose que $\Chi(S) \le 0$ et que
$\Chi(V,S) \le -\Chi(S)$.

\alinea{\textbf{a)}} 
Il existe sur~$V$ des structures de contact transversales aux fibres et
d'enroulement~$-n$ si et seulement si $n=1$ ou si $n \Chi(V,S) = -\Chi(S)$ et
$n>0$.

\alinea{\textbf{b)}}
Si $\Chi(V,S) \ne -\Chi(S)$, les structures de contact universellement tendues
et d'enroulement~$-1$ sur~$V$ forment une seule classe d'isotopie.

\alinea{\textbf{c)}}
Si $n \Chi(V,S) = -\Chi(S)$, $n>0$, les structures de contact universellement
tendues et d'enroulement~$-n$ forment un nombre fini de classes de conjugaison
\'egal au nombre de diviseurs de~$n$. De plus, chaque classe de conjugaison
contient une infinit\'e de classes d'isotopie.
\end{theoreme}

Avant de d\'emontrer ce th\'eor\`eme dans les sections C---F, on observe que la
partie~c) pr\'esente peu d'int\'er\^et lorsque $\Chi(S) = 0$. En effet, comme
$n \Chi(V,S) = -\Chi(S)$, la vari\'et\'e~$V$ est un tore et l'enroulement n'est
pas invariant par conjugaison. Dans ce cas, le bon invariant de conjugaison est
la torsion~\cite{Gi:bifurcations}. D'autre part, le cas laiss\'e de c\^ot\'e par
le th\'eor\`eme est en fait beaucoup plus simple :

\begin{proposition} \label{t:sphere}
Soit $V$ une vari\'et\'e connexe et orient\'ee, fibr\'ee en cercles au-dessus
d'une sph\`ere~$S$. On suppose que $\Chi(V,S) < 0$. \`A isotopie pr\`es, il y a
sur~$V$ une seule structure de contact directe et transversale aux fibres. Son
enroulement vaut~$-2$ si $V \simeq \S^3$ et~$-1$ sinon.
\end{proposition}

\begin{proof}
Ce r\'esultat est inclus dans le th\'eor\`eme~1.1 de \cite{Gi:bifurcations} qui
classifie les structures de contact tendues sur les espaces lenticulaires. Comme
il n'en constitue qu'une toute petite partie, on indique bri\`evement sa preuve.
Soit $\xi_0$ et $\xi_1$ deux structures de contact sur~$V$ qui sont directes et
transversales aux fibres. Quitte \`a d\'eformer l'une d'elles par une isotopie
qui la laisse transversale aux fibres, on peut supposer que $\xi_0$ et~$\xi_1$
co\"{\i}ncident au-dessus de la r\'eunion disjointe~$Q$ de deux petits disques
dans~$S$. On param\`etre alors $V \setminus \pi^{-1}(Q)$ par $\T^2 \times [0,1]$
de telle sorte que $\pi$ soit la projection $\T^2 \times [0,1] \to \S^1 \times
[0,1]$, \ $(x,y,t) \mapsto (y,t)$. Chaque structure $\xi_i$ imprime ainsi sur
chaque tore $\T^2 \times \{t\}$ un feuilletage caract\'eristique dont toutes les
feuilles coupent transversalement les fibres de~$\pi$. Dans la terminologie de
\cite{Gi:bifurcations}, $\xi_0$ et $\xi_1$ sont des structures rotatives sur
$\T^2 \times [0,1]$ qui co\"\i ncident pr\`es du bord et ont la m\^eme amplitude
(non nulle). Le th\'eor\`eme~3.3 de \cite{Gi:bifurcations} montre qu'elles sont
alors isotopes relativement au bord.
\end{proof}

\subsection{Fibrations legendriennes et rev\^etements fibr\'es}

On \'etudie ici les structures de contact sur~$V$ qui sont tangentes aux fibres,
\emph{i.e.} pour lesquelles la projection $\pi \from V \to S$ est une fibration
legendrienne. Comme la proposition~\ref{p:perturbation} le laisse entrevoir, ces
structures jouent un r\^ole important dans l'\'etude des structures de contact
transversales aux fibres. On note cependant qu'une structure de contact tangente
aux fibres n'est pas n\'ecessairement orientable le long des fibres. L'exemple
type est la structure de contact canonique~$\xi_S$ sur le fibr\'e $\P(TS)$ des
droites non orient\'ees tangentes \`a~$S$.

\begin{proposition} \label{p:legendre}
Soit $V$ une vari\'et\'e connexe et orient\'ee, fibr\'ee en cercles au-dessus
d'une surface quelconque~$S$. L'application qui, \`a tout rev\^etement fibr\'e
$\rho \from V \to \P(TS)$, associe la structure de contact $\xi = \rho^*\xi_S$
est une bijection de l'espace des rev\^etements fibr\'es et orient\'es $V \to
\P(TS)$ dans l'espace des structures de contact tangentes aux fibres sur~$V$. En
outre, l'enroulement de $\rho^*\xi_S$ autour des fibres vaut~$-d/2$ o\`u $d$ est
le degr\'e du rev\^etement $\rho$.
\end{proposition}

\begin{proof}
Si $\rho \from V \to \P(TS)$ est un rev\^etement fibr\'e orient\'e, la structure
de contact $\xi = \rho^*\xi_S$ est tangente aux fibres de $V$ (et directe). En
outre, comme l'enroulement de~$\xi_S$ autour des fibres vaut $-1/2$, celui de
$\xi$ vaut $-d/2$ o\`u $d$ est le degr\'e de~$\rho$. D'autre part, tout champ de
plans $\xi$ tangent aux fibres de~$V$ d\'efinit une application fibr\'ee $V \to
\P(TS)$ : l'image d'un point~$p$ est simplement la projection sur~$S$ du plan
$\xi(p)$. Cette application est un rev\^etement (orient\'e) si et seulement si
le champ~$\xi$ est une structure de contact (directe) et cette structure est,
par construction, le rappel de~$\xi_S$.
\end{proof}

\begin{corollaire} \label{c:legendre}
Soit $V$ une vari\'et\'e connexe et orient\'ee fibr\'ee en cercles au-dessus
d'une surface close~$S$. Pour que $V$ porte une structure de contact tangente
aux fibres \up(et directe\up), il faut et il suffit qu'il existe un entier~$d>0$
tel que $d \Chi(V,S) = -2\Chi(S)$.
\end{corollaire}

\begin{proof}
Le nombre d'Euler du fibr\'e $\P(TS) \to S$, pour l'orientation induite par la
structure de contact~$\xi_S$, est $-2\Chi(S)$. Ainsi, la relation $d \Chi(V,S) =
-2\Chi(S)$ traduit simplement l'existence d'un rev\^etement fibr\'e et orient\'e
de $V$ sur $\P(TS)$ \`a $d$~feuillets.
\end{proof}

\smallskip

\begin{remarque}
Lorsque la surface~$S$ est close, l'existence sur~$V$ d'un simple champ de plans
tangent aux fibres exige en fait d\'ej\`a que le rapport $-2\Chi(S)/\Chi(V,S)$
soit un entier, \'eventuellement n\'egatif ou nul. En effet, ce champ de plans
d\'efinit une application fibr\'ee $V \to \P(TS)$ dont la restriction \`a chaque
fibre est de degr\'e constant~$d$. D'autre part, un feuilletage tangent aux
fibres ne peut \^etre que l'image inverse d'un feuilletage de~$S$ et n'existe
donc que si $\Chi(S) = 0$.
\end{remarque}

\begin{corollaire} \label{c:whitney}
Soit $V$ une vari\'et\'e connexe et orient\'ee, fibr\'ee en cercles au-dessus
d'une surface compacte~$S$ ayant un bord connexe non vide. On note $\wh S
\subset V$ l'image d'une section, $L \subset \partial V$ une fibre et on oriente
$\partial \wh S$ et~$L$ de telle sorte que leur intersection soit positive sur
$\partial V$.

Pour toute structure de contact~$\xi$ tangente aux fibres de~$V$, la
multi-courbe des singularit\'es de $\xi\, \partial V$ a pour classe d'homologie
$$ \pm2 \bigl( \w(\xi) \, [\partial \wh S] + \Chi(S) \, [L] \bigr) \in
   H_1 (\partial V;\Z), $$
du moins si toutes ses composantes sont orient\'ees dans le m\^eme sens.
\end{corollaire}

\begin{proof}
Si $V = \P(TS)$ et si $\xi$ est la structure de contact canonique, la courbe des
singularit\'es de $\xi\, \partial V$ n'est autre que le rel\`evement legendrien
de $\partial S$ dans~$V$. Comme $\w(\xi) = -1/2$, la formule dit simplement que,
par rapport \`a n'importe quel champ de droites d\'efini sur~$S$, la tangente au
bord~$\partial S$ fait $-2\Chi(S)$~tours sur la fibre\footnote
{Pour trouver le bon signe, noter que l'orientation de contact sur $\P(TS)$ en
un point $(q,\delta)$, $\delta \subset T_qS$, n'est pas la juxtaposition d'une
orientation de $T_qS$ et de l'orientation induite sur la droite $\P(T_qS)$, mais
l'inverse.}.
Le cas g\'en\'eral s'obtient en passant \`a un rev\^etement fibr\'e de degr\'e
$-2\w(\xi)$.
\end{proof}

\smallskip

Pour terminer cette section, on v\'erifie que l'enroulement d'une structure de
contact tangente aux fibres est bien ce qu'on attend :

\begin{lemme} \label{l:bennequin}
Soit $V$ une vari\'et\'e connexe et orient\'ee, fibr\'ee en cercles au-dessus
d'une surface quelconque~$S$, et soit $\xi$ une structure de contact tangente
aux fibres. L'enroulement de~$\xi$ est \'egal \`a son enroulement autour des
fibres.
\end{lemme}

\begin{proof}
Quitte \`a se placer au-dessus du rev\^etement universel de~$S$, on suppose que
$S$ est un plan ou une sph\`ere. Si $S \simeq \S^2$, il suffit de traiter le cas
o\`u $V \simeq \S^3$. Dans ce cas, l'enroulement de~$\xi$ autour des fibres est
\'egal \`a~$-2$. Par ailleurs, comme on l'a not\'e au d\'ebut de la section~2.B,
l'enroulement de $\xi$ autour de toute courbe legendrienne $L$ isotope \`a la
fibre dans $\S^3$ ---~donc non nou\'ee~--- vaut $\tb(L)-1$. L'identit\'e voulue
r\'esulte donc de l'in\'egalit\'e de Bennequin $\tb(L) \le -1$.

Si $S \simeq \R^2$, il existe un diff\'eomorphisme $V \to \R^2 \times \S^1$ qui
conjugue la fibration~$\pi$ \`a la projection sur~$\R^2$ et envoie~$\xi$ sur la
structure de contact d'\'equation
$$ \cos(n\theta)\, dx - \sin(n\theta)\, dy = 0, \qquad
   (x,y,\theta) \in \R^2 \times \S^1, $$
o\`u $n$~est l'enroulement de~$\xi$ autour des fibres. Dans $\R^3$ muni de sa
structure de contact ordinaire, on consid\`ere alors une courbe legendrienne non
nou\'ee~$L$ dont l'invariant de Thurston-Bennequin vaut~$-1$. Un avatar du
th\'eor\`eme de Darboux donne un plongement de $(V,\xi)$ dans $\R^3$ qui envoie
sur~$L$ une fibre de~$\pi$. Toute courbe legendrienne dans $(V,\xi)$ isotope \`a
la fibre et d'enroulement~$-m>-n$ a alors pour image une courbe legendrienne non
nou\'ee dont l'invariant de Thurston-Bennequin vaut $n-m-1$. L'in\'egalit\'e de
Bennequin permet \`a nouveau de conclure.
\end{proof}

\subsection{Mise en position tangentielle}

On d\'emontre ici la partie~a) du th\'eor\`eme~\ref{t:ghys}. On observe d'abord
que la structure de contact canonique sur la vari\'et\'e $\S(TS)$ des droites
orient\'ees tangentes \`a~$S$ a pour enroulement~$-1$ (lemme~\ref{l:bennequin}).
Par suite, comme $\Chi(V,S) \le -\Chi(S)=\Chi(\S(TS),S)$, la chirurgie d\'ecrite
dans le lemme~\ref{l:chirurgie} permet de produire sur~$V$ une structure de
contact transversale aux fibres et d'enroulement~$-1$. D'autre part, lorsque $n
\Chi(V,S) = -\Chi(S)$ pour un certain entier~$n$, la vari\'et\'e~$V$ admet un
rev\^etement fibr\'e \`a $n$~feuillets au-dessus de $\S(TS)$. Ainsi, $V$~porte
une structure de contact tangente aux fibres, d'enroulement~$-n$ en vertu du
lemme~\ref{l:bennequin}, que la proposition~\ref{p:perturbation} permet de
rendre transversale aux fibres. Il reste donc \`a d\'emontrer la proposition
suivante :

\begin{proposition} \label{p:ghys}
Soit $V$ une vari\'et\'e connexe et orient\'ee, fibr\'ee en cercles au-dessus
d'une surface close orientable~$S$, et soit~$\xi$ une structure de contact
transversale aux fibres et d'enroulement $\w(\xi) = -n$, $n \in \N$. Si $n
\Chi(V,S) = -\Chi(S)$, alors $\xi$ est isotope \`a une structure de contact
tangente aux fibres. Sinon, $n$~vaut~$1$.
\end{proposition}

\begin{proof}
On note $D \subset S$ un disque ferm\'e, $R$ la surface $S \setminus \Int D$ et
$W$ le tore plein $\pi^{-1}(D)$. On param\`etre~$W$ par $\D^2 \times \S^1$ de
telle sorte que la fibration $\pi\res W$ soit la projection sur $\D^2$. En vertu
de la proposition~\ref{p:calculs} et du lemme~\ref{l:roussarie}, $\xi$ est
isotope \`a une structure de contact $\xi'$ ayant les propri\'et\'es suivantes :
\begin{itemize}
\item
au-dessus de~$R$, les fibres sont tangentes \`a~$\xi'$ et ont pour enroulement
$\w(\xi') = \w(\xi)$ ;
\item
au-dessus de $\Int D$, les fibres sont transversales \`a~$\xi'$.
\end{itemize}
D'apr\`es le lemme~\ref{l:deux}, les singularit\'es du feuilletage $\xi' \,
\partial W$ forment deux cercles. De plus, vu la d\'efinition du nombre d'Euler
$\Chi(V,S)$, le corollaire~\ref{c:whitney} montre que la pente de ces cercles
sur $\partial W = \partial \D^2 \times \S^1$ ---~ou de leurs classes d'homologie
dans $H_1(\partial W;\R) = \R^2$~--- vaut
$$ \mu = \frac 1 n \bigl( n \Chi(V,S) + \Chi(S) - 1 \bigr) \,. $$

\begin{assertion}
Si $n>1$, alors $n \Chi(V,S) = - \Chi(S)$.
\end{assertion}

\begin{proof}[Preuve]
Soit $T_a$, $0<a<1$, le tore $a\S^1 \times \S^1 \subset W = \D^2 \times \S^1$.
Le feuilletage caract\'eristique $\xi' T_a$ de chaque tore $T_a$ est transversal
aux fibres et est donc d\'ecrit par l'application~$\phi_a$ de premier retour sur
une fibre. Quand $a$ varie de~$0$ \`a~$1$, le nombre de translation $\mu_a$ de~$
\phi_a$ d\'ecro\^\i t contin\^ument de $0$ \`a~$\mu$ (proposition~\ref{p:wood}).
En particulier, $\mu$ est n\'egatif ou nul. Par suite, comme $n = -\w(\xi) > 0$
(proposition~\ref{p:calculs}), l'entier $n \Chi(V,S) + \Chi(S)$ est n\'egatif
ou nul.

Si $n>1$ et $n\Chi(V,S) + \Chi(S) < 0$, la pente~$\mu$ est strictement major\'ee
par $-1/(n-1)$. D'autre part, une petite perturbation de~$\xi'$ dans $W$ permet
d'imposer \`a la famille $\phi_a$ n'importe quelle propri\'et\'e g\'en\'erique.
On peut ainsi supposer que, pour une valeur~$a$, le diff\'eomorphisme it\'er\'e
$\phi_a^{n-1}$ a pour nombre de translation $(n-1) \mu_a = -1$ et ne poss\`ede
que deux points fixes, lesquels sont hyperboliques. Le tore~$T_a$ correspondant
est alors convexe dans $(V,\xi')$ et on note $\Gamma \subset T_a$ une bi-courbe
transversale aux fibres qui scinde le feuilletage $\xi' T_a$. On choisit ensuite
sur~$T_a$ un feuilletage singulier~$\sigma$ scind\'e par~$\Gamma$ et pour lequel
chaque fibre~$L$ de $\pi \res {T_a}$ est satur\'ee (voir la d\'emonstration du
lemme~\ref{l:roussarie}). Le lemme \ref{l:ctc} fournit un plongement $\phi \from
T_a \to W$ isotope \`a l'inclusion et dont l'image $T = \phi(T_a)$ a pour
feuilletage caract\'eristique $\xi' T = \phi_*\sigma$. Ainsi, chaque courbe
$\phi(L)$ est isotope \`a la fibre et legendrienne. De plus, comme $\Gamma$ a
deux composantes, l'enroulement de $\phi(L)$ vaut $1-n > -n = \w(\xi)$, ce qui
est absurde.
\end{proof}

\smallskip

On suppose d\'esormais que $n \Chi(V,S) = -\Chi(S)$ de sorte que $\mu = -1/n$.
Pour finir la d\'emonstration de la proposition~\ref{p:ghys}, il reste \`a
montrer que~$\xi$ ---~ou $\xi'$~--- est isotope \`a une structure de contact
partout tangente aux fibres. Pour cela, on note d'abord que, d'apr\`es le lemme
\ref{l:univ}, la restriction de $\xi'$ \`a~$W$ est universellement tendue. D\`es
lors, le th\'eor\`eme~1.6 de \cite{Gi:bifurcations} assure que toute structure 
de contact universellement tendue sur~$W$ qui imprime le m\^eme feuilletage que
$\xi'$ sur $\partial W$ est isotope \`a $\xi'$ relativement au bord. En effet,
la condition $\mu = -1/n$ garantit  qu'il existe dans~$W$, \`a isotopie relative
au bord pr\`es, un seul anneau s'appuyant sur les cercles de singularit\'es de
$\xi'\, \partial W$ (voir le d\'ebut de la d\'emonstration de la proposition
\ref{p:bfs}). Or la structure d'\'equation 
$$ \cos(n\theta) \, dx - \sin(n\theta) \, dy = 0, \qquad
   (x,y,\theta) \in \D^2 \times \S^1, $$
d\'efinit sur $W = \D^2 \times \S^1$ une structure de contact universellement
tendue, tangente \`a toutes les fibres et d'enroulement~$-n$ autour de chacune.
\end{proof}

\subsection{Structures d'enroulement l\^ache}

On d\'emontre ici la partie~b) du th\'eor\`eme~\ref{t:ghys}. On suppose donc que
$\Chi(V,S) \ne -\Chi(S)$ et on consid\`ere sur~$V$ deux structures de contact
$\xi_0$ et~$\xi_1$ universellement tendues et d'enroulement~$-1$. On note~$D$ un
disque ferm\'e dans~$S$ et on regarde la surface $R = S \setminus \Int D$ comme
le voisinage r\'egulier d'un bouquet de cercles $K = \bigvee_{i=1}^{2g}K_i$ o\`u
$g$ est le genre de~$S$. On param\`etre en outre le tore plein $W = \pi^{-1}(D)$
par $\D^2 \times \S^1$ de telle sorte que la fibration $\pi \res W$ soit la
projection sur~$\D^2$.

Compte tenu des lemmes~\ref{l:roussarie}, \ref{l:univ} et de la proposition~\ref
{p:bfs}, on suppose que chaque structure de contact~$\xi_i$, $i \in \{0,1\}$,
satisfait aux conditions suivantes :
\begin{itemize}
\item
au-dessus de~$R$, les fibres sont tangentes \`a~$\xi_i$ et d'enroulement~$-1$ ;
\item
au-dessus de~$\Int D$, les fibres sont transversales \`a~$\xi_i$.
\end{itemize}
Les singularit\'es du feuilletage $\xi_i\,\partial W$ forment alors deux cercles
dont la classe d'homologie vaut $\pm (1,\Chi(V,S)+\Chi(S)-1)$, o\`u $\Chi(V,S) +
\Chi(S) - 1 \le -2$.

Soit $Q \subset R$ un voisinage compact du sommet de $K$ dont l'intersection
avec le bouquet~$K$ est connexe. On lisse chaque lacet~$K_i$ dans $Q$ en une
courbe~$K'_i$ et on param\`etre chaque tore $F_i = \pi^{-1}(K'_i)$ par $\T^2$ de
telle sorte que les fibres soient les cercles $\{*\} \times \S^1$ et que les
(deux) courbes de singularit\'es du feuilletage $\xi_0 F_i$ aient pour classe
d'homologie $\pm(1,0)$. Les courbes de singularit\'es du feuilletage $\xi_1 F_i$
ont alors une classe du type $\pm(1,n_i)$, $n_i \in \Z$, et, quitte \`a composer
le param\'etrage de $F_i$ par la transformation $(x_1,x_2) \in \T^2 \mapsto
(-x_1,x_2)$, on prend $n_i\ge0$.

\begin{lemme} \label{l:lache}
Il existe une structure de contact isotope \`a~$\xi_1$ qui co\"{\i}ncide avec~$
\xi_0$ au-dessus de~$R$.
\end{lemme}

\begin{proof}
Si les entiers~$n_i$ sont tous nuls, une isotopie fibr\'ee am\`ene~$\xi_1$ \`a
co\"{\i}ncider avec $\xi_0$ au-dessus de~$R$. On suppose ci-dessous $n_i\ne0$ et
on d\'eforme~$\xi_1$, par une isotopie relative \`a $\bigvee_{j \ne i} F_j$, en
une structure de contact $\xi_2$ dans laquelle $F_i$ est convexe et poss\`ede un
feuilletage caract\'eristique $\xi_2 F_i$ scind\'e par deux courbes, de classe
$\pm(1,n_i-1)$. Comme dans la d\'emonstration du lemme~\ref{l:roussarie}, les
lemmes sur les surfaces convexes permettent ensuite de modifier~$\xi_2$, par une
isotopie toujours relative \`a $\bigvee_{j \ne i} F_j$ au cours de laquelle
$F_i$ reste convexe, en une structure de contact~$\xi_3$ tangente aux fibres de
$\pi$ sur un voisinage de~$F_i$. En renouvelant l'op\'eration, on annule au fur
et \`a mesure les coefficients~$n_i$.

\smallskip

Soit $J_s$, $s \in [0,1]$, une famille lisse d'arcs plong\'es dans $S \setminus
(\Int Q \cup \bigvee_{j \ne i} K_j)$ et ayant les propri\'et\'es suivantes :
\begin{itemize}
\item
$J_0 = K_i \setminus \Int Q$ et $J_s$, pour tout $s \in [0,1]$, a un contact
d'ordre infini avec~$J_0$ en ses extr\'emit\'es ;
\item
les arcs $J_s$, $s \in [0,1]$, sont d'int\'erieurs disjoints, leur union couvre
un voisinage de~$D$ et la trace de chacun d'eux sur~$\Int D$ est connexe.
\item
les anneaux $B_s = \pi^{-1}(J_s)$, pour $s>0$, sont du c\^ot\'e positif de $F_i$
---~la coorientation provient du param\'etrage $\T^2 \to F_i$.
\end{itemize}
Pour tout $s \in [0,1]$, on pose $F_{i,s} = (F_i \setminus B_0) \cup B_s$. Le
tore $F_{i,1}$ est contenu dans $\pi^{-1}(R)$ et les cercles de singularit\'es
du feuilletage $\xi_1 F_{i,1}$ ont pour classe $\pm(1,n_i+\Chi(V,S)+\Chi(S))$.
En fait, chaque feuilletage $\xi_1 F_{i,s}$ est tangent aux fibres au-dessus de
$J_s \setminus D$ (avec deux singularit\'es par fibre) et transversal aux fibres
au-dessus de $J_s \cap \Int D$. Les singularit\'es de $\xi_1 F_{i,s}$ forment
donc deux courbes $C_s^\pm$ qui sont ferm\'ees si $J_s$ \'evite $\Int D$ (par
exemple pour $s$ proche de $0$ et $1$) mais sont des arcs sinon. Sauf pour un
nombre fini de valeurs de~$s$ (instants de bifurcation), les extr\'emit\'es de
l'arc $C_s^+$ sont reli\'ees par une feuille r\'eguli\`ere. La courbe ferm\'ee
$\ol C_s^+$, union de cette feuille et de $C_s^+$, a une classe d'homologie qui,
lorsque $s$ varie de $0$ \`a~$1$, prend successivement les valeurs $\pm(1,n_i)$,
$\pm(1,n_i-1),\dotsc$, $\pm(1,n_i+\Chi(V,S)+\Chi(S))$ \ (la condition de contact
impose la d\'ecroissance de la pente, comme dans la proposition~\ref{p:wood}).

On choisit d\'esormais pour $s$ un instant o\`u la classe de $\ol C_s^+$ vaut
$\pm(1,n_i-1)$. Le feuilletage $\xi_1 F_{i,s}$ est clairement scind\'e par le
bord d'un voisinage annulaire de~$\ol C_s^+$ donc $F_{i,s}$ est convexe. Soit
$\phi_t$, $t \in [0,s]$, une isotopie de $V$ qui prolonge l'isotopie $F_{i,t}$
sans bouger les points de $\pi^{-1}(Q) \cup \bigvee_{j\ne i}F_j$. Les structures
de contact $\phi_t^*\xi_1$ donnent la d\'eformation voulue entre~$\xi_1$ et
$\xi_2 = \phi_s^*\xi_1$.
\end{proof}

\smallskip

Fort du lemme~\ref{l:lache}, on suppose dor\'enavant que $\xi_1$ co\"{\i}ncide
avec~$\xi_0$ au-dessus de~$R$ et on pose $\sigma = \xi_i\, \partial W$. On note
que, d'apr\`es le lemme~\ref{l:univ}, les restrictions de $\xi_0$ et~$\xi_1$
\`a~$W$ sont universellement tendues. Or, d'apr\`es le th\'eor\`eme~1.6 de \cite
{Gi:bifurcations}, les structures de contact universellement tendues sur~$W$ qui
impriment~$\sigma$ sur $\partial W$ forment deux classes d'isotopie relative au
bord. Mieux, $\xi_0 \res W$ et $\xi_1 \res W$ sont dans la m\^eme classe si et
seulement si les anneaux respectifs $A_0$ et~$A_1$ que leur attribue le lemme
\ref{l:bfs} ---~anneaux qui s'appuient sur les cercles de singularit\'es de
$\sigma$~--- sont isotopes relativement \`a leur bord (cf. d\'emonstration de la
proposition~\ref{p:bfs}). Si $A_1$ n'est pas isotope \`a~$A_0$ relativement
\`a son bord, on am\`ene $A_1$ sur~$A_0$ par une isotopie fibr\'ee de~$W$ qui
permute les cercles de singularit\'es de~$\sigma$. On prolonge ci-dessous cette
isotopie en une isotopie fibr\'ee~$\phi_t$ de $V$ dont le stade final~$\phi_1$
pr\'eserve $\xi_1 \res {\pi^{-1}(R)}$. Du coup, $\xi_0$ est isotope $(\phi_1)_*
\xi_1$ relativement \`a $\pi^{-1}(R)$, donc $\xi_0$ et $\xi_1$ sont isotopes.

Comme $R$ est orientable, un avatar de la proposition~\ref{p:legendre} donne un
diff\'eomorphisme fibr\'e de $\pi^{-1}(R)$ sur $\S(TR)$ qui envoie $\xi_1$ sur
la structure de contact canonique~$\xi_R$. Moyennant le choix d'une structure
conforme sur $R$, on note $\psi_t$ l'isotopie fibr\'ee de $\S(TR)$ qui tourne
les droites d'un angle $\pi t$, $t \in [0,1]$. Le diff\'eomorphisme~$\psi_1$
pr\'eserve~$\xi_R$ et permute les courbes de singularit\'es du feuilletage
$\xi_R\, \partial \S(TR)$. Par suite, l'isotopie $\psi_t$, transport\'ee sur
$\pi^{-1}(R)$, fournit le prolongement voulu.
\qed

\subsection{Structures d'enroulement serr\'e}

On d\'emontre ici la partie~c) du th\'eor\`eme~\ref{t:ghys}, au calcul pr\`es du
nombre exact des classes de conjugaison qu'on effectue dans la proposition~\ref
{p:otal}. On suppose donc que $n \Chi(V,S) = -\Chi(S)$, o\`u $n>0$ est entier,
et on s'int\'eresse sur~$V$ aux structures de contact universellement tendues et
d'enroulement~$-n$. Toute structure de ce type est orientable le long des fibres
(car $n$ est entier) et est isotope, d'apr\`es la proposition~\ref{p:ghys}, \`a
une structure de contact tangente aux fibres. D'autre part, l'avatar orient\'e
de la proposition~\ref{p:legendre} assure que toute structure de contact~$\xi$
tangente aux fibres et orientable le long des fibres s'\'ecrit $\rho^*\xi_S$,
o\`u $\rho$ est un rev\^etement fibr\'e $V \to \S(TS)$ et $\xi_S$ la structure
de contact canonique sur $\S(TS)$.

On dira que deux rev\^etements fibr\'es $\rho_0, \rho_1 \from V \to \S(TS)$ sont
\emph{isomorphes au-dessus de~$S$} s'il existe des diff\'eomorphismes $\phi$ de
$V$ et $\ol\phi$ de~$\S(TS)$, fibr\'es au-dessus d'un diff\'eomorphisme de~$S$,
qui rendent commutatif le diagramme
$$ \begin{CD}
   V    @>\phi>>    V    \\
@V{\rho_0}VV @VV{\rho_1}V \\
\S(TS) @>>{\ol\phi}> \S(TS) \rlap{\,.}
\end{CD} $$

\smallskip

Comme dans la section~D, on prend sur $S$ un bouquet de cercles $K
= \bigvee_{i=1}^{2g} K_i$ ayant pour compl\'ementaire un disque et on lisse les
lacets~$K_i$ en des courbes~$K'_i$. On param\`etre de nouveau chaque tore $F_i =
\pi^{-1}(K'_i)$ par $\T^2$ de telle sorte que les fibres soient les cercles
$\{*\} \times \S^1$. Pour toute structure de contact~$\xi$ sur~$V$ tangente
aux fibres et d'enroulement~$-n$, la courbe des singularit\'es du feuilletage
$\xi F_i$ ---~toutes composantes orient\'ees dans le m\^eme sens~--- a alors une
classe d'homologie qui s'\'ecrit
$$ \pm2 \bigl( n, m_i(\xi) \bigr) \in H_1(F_i;\Z) \cong \Z^2 \,. $$
En outre, \'etant donn\'e des entiers $m_1, \dots, m_{2g} \in \Z$, on construit
sans peine une structure de contact $\xi$ tangente aux fibres et d'enroulement
$-n$ pour laquelle $m_i = m_i(\xi)$, $1 \le i \le 2g$. Lorsque $\Chi(S) < 0$, la
partie~c) du th\'eor\`eme~\ref{t:ghys} d\'ecoule donc du lemme suivant et de la
proposition~\ref{p:otal}) :

\begin{lemme} \label{l:serre}
Soit $\xi_0$ et $\xi_1$ deux structures de contact tangentes aux fibres et
d'enroulement~$-n$. 

\alinea{\textbf{a)}}
Les structures $\xi_0$ et $\xi_1$ sont isotopes si et seulement si les entiers
$m_i(\xi_0)$ et $m_i(\xi_1)$ sont \'egaux pour $1 \le i \le 2g$.

\alinea{\textbf{b)}}
Lorsque $\Chi(S) < 0$, les structures $\xi_0$ et $\xi_1$ sont conjugu\'ees si et
seulement si les rev\^etements fibr\'es associ\'es, de $V$ sur $\S(TS)$, sont
isomorphes au-dessus de~$S$.
\end{lemme}

\begin{proof}
\Alinea{a)}
Si $m_i(\xi_0) = m_i(\xi_1)$ pour $1 \le i \le 2g$, une isotopie fibr\'ee permet
clairement d'amener~$\xi_0$ sur~$\xi_1$. Par ailleurs, si on note~$q$ le sommet
du bouquet~$K$, chaque lacet~$K_i$ engendre dans $\pi_1(S,q)$ un groupe cyclique
infini. Il lui est ainsi associ\'e un rev\^etement $\rho_i \from \wt S_i \to S$
et on param\`etre $\wt S_i$ par $\R \times \S^1$ de telle sorte que $\rho(\{0\}
\times \S^1)$ soit la courbe $K'_i$. Les rel\`evements respectifs $\wt\xi_0$ et
$\wt\xi_1$ de $\xi_0$ et~$\xi_1$ sur $\wt V = \rho^*V$ sont des structures de
contact tangentes aux fibres de la fibration $\wt V \to \wt S$. De plus, on peut
param\'etrer $\wt V$ par $\R \times \T^2$ de telle sorte que les conditions
suivantes soient remplies :
\begin{itemize}
\item
le plongement compos\'e $\T^2 = \{0\} \times \T^2 \to \wt V \to V$ a pour image
$F_i$ et co\"{\i}ncide avec le param\'etrage donn\'e de $F_i$ ;
\item
la fibration $\wt V \to \wt S$ est la projection sur $\R \times \S^1$.
\end{itemize}
Du coup, tous les tores $T_a = \{a\} \times \T^2$, $a \in \R$, sont convexes et
la courbe qui scinde leur feuilletage $\wt \xi_0 T_a$ (resp. $\wt\xi_1 T_a$) a
pour classe $\pm (1,m_i(\xi_0))$ (resp. $\pm (1,m_i(\xi_1))$). Si $m_i(\xi_0)$
est diff\'erent de $m_i(\xi_1)$, le lemme~4.7 de \cite{Gi:bifurcations} (qui est
un cas particulier de l'in\'egalit\'e de Bennequin semi-locale, cf. proposition
\ref{p:bennequin}) montre que $\wt\xi_0$ n'est pas isotope \`a~$\wt\xi_1$. Par
suite, $\xi_0$ et $\xi_1$ ne sont pas isotopes.

\smallskip

\alinea{b)}
Soit $\rho_0$ et $\rho_1$ les rev\^etements fibr\'es $V \to \S(TS)$ associ\'es
respectivement \`a $\xi_0$ et $\xi_1$. Si $\rho_0$ et $\rho_1$ sont isomorphes
au-dessus de~$S$, tout diff\'eomorphisme fibr\'e de~$V$ qui les conjugue envoie
en m\^eme temps $\xi_0$ sur $\xi_1$. On suppose donc maintenant que $\xi_0$ et
$\xi_1$ sont conjugu\'ees par un diff\'eomorphisme $\phi_0$ de $V$. D'apr\`es
\cite{Wa}, $\phi_0$ est isotope \`a un diff\'eomorphisme fibr\'e $\phi$. La
structure de contact~$\phi^*\xi_1$ est alors tangente aux fibres et isotope \`a
$\xi_0$. Or il ressort imm\'ediatement du~a) que, si deux structures de contact
tangentes aux fibres et de m\^eme enroulement sont isotopes, elles le sont par
une isotopie fibr\'ee. \`A l'instant final, cette isotopie conjugue $\rho_0$
\`a $\phi^*\rho_1$, donc $\rho_0$ et $\rho_1$ sont isomorphes au-dessus de~$S$.
\end{proof}

\smallskip

Il reste \`a regarder le cas (peu int\'eressant) o\`u $\Chi(S) = 0$. La relation
$n \Chi(V,S) = -\Chi(S)$ force alors $V$ \`a \^etre un tore. D'apr\`es
\cite{Gi:bifurcations} (voir aussi \cite{Ka}), toute structure de contact
(universellement) tendue sur $\T^3$ est conjugu\'ee \`a une structure de contact
$\zeta_m$ d'\'equation
$$ \cos(m\theta)\, dx_1 - \sin(m\theta)\, dx_2 = 0, \qquad
   m>0, \ \ (x_1,x_2,\theta) \in \T^3 \,. $$
Pour la projection $(x_1,x_2,\theta) \mapsto (x_1,x_2)$, la structure~$\zeta_m$
est tangente aux fibres et d'enroulement~$-m$. Mais, pour chaque entier $d>0$,
on peut aussi trouver une fibration $\T^3 \to \T^2$ pour laquelle $\zeta_m$ est
d'enroulement~$-dm$.
\qed

\subsection{D\'enombrement des rev\^etements fibr\'es}

On classifie ici les rev\^etements fibr\'es de $\S(TS)$ \`a isomorphisme pr\`es
au-dessus de~$S$, ce qui compl\`ete la d\'emonstration du th\'eor\`eme~\ref
{t:ghys}-c.

\begin{proposition} \label{p:otal}
Soit $S$ une surface close orientable, de caract\'eristique d'Euler n\'ega\-tive
ou nulle et soit $n$~un entier positif.

\alinea{\textbf{a)}}
\`A isomorphisme pr\`es au-dessus de~$S$, les rev\^etements fibr\'es de $\S(TS)$
\`a $n$~feuillets sont class\'es par le quotient $H^1(S;\Z/n\Z) \,\big/ \Aut
\pi_1(S)$. 

\alinea{\textbf{b)}}
\`A composition pr\`es par les automorphismes de $\pi_1(S)$, les morphismes de
$\pi_1(S)$ dans $\Z/n\Z$ sont class\'es par leur image. Le cardinal du quotient
$H^1(S;\Z/n\Z) \,\big/ \Aut \pi_1(S)$ est donc \'egal au nombre de diviseurs de
$n$.
\end{proposition}

\begin{proof}
\Alinea{a)}
On pose $V_0 = \S(TS)$ et on choisit dans~$V_0$ un point de r\'ef\'erence. Comme
$\Chi(S) \le 0$, le groupe fondamental $\pi_1(V_0)$ est une extension centrale
de $\pi_1(S)$ par $\Z \cong \pi_1(\S^1)$. La suite exacte courte d'homotopie $\Z
\to \pi_1(V_0) \to \pi_1(S)$ induit ainsi une suite exacte courte d'homologie
$\Z/\Chi(S)\Z \to H_1(V_0;\Z) \to H_1(S)$.

D'autre part, la bijection entre les (classes d'\'equivalence de) rev\^etements
de~$V_0$ et les (classes de conjugaison de) sous-groupes de $\pi_1(V_0)$ associe
aux rev\^etements fibr\'es de degr\'e~$n$ les sous-groupes d'indice~$n$ dont la
trace sur le sous-groupe central $\Z$ est $n\Z$. Ces sous-groupes sont normaux
et sont donc les noyaux des morphismes $\pi_1(V_0) \to \Z/n\Z$ dont la
restriction
\`a $\Z$ est la projection canonique. Comme de tels morphismes transitent par
$H_1(V_0;\Z)$, ils existent d\`es que $n$ divise $\Chi(S)$ et forment alors un
espace principal homog\`ene du groupe $H^1(S;\Z/n\Z)$. En outre, l'action des
diff\'eomorphismes fibr\'es de $V_0$ se r\'eduit sur $H^1(S;\Z/n\Z)$ \`a
l'action des diff\'eomorphismes de $S$, \emph{i.e.} des automorphismes
(ext\'erieurs) de $\pi_1(S)$.

\smallskip

\alinea{b)}
Soit $f_0$ et $f_1$ des morphismes de $\pi_1(S)$ dans $\Z/n\Z$. Comme $\pi_1(S)$
est sans torsion, $f_0$ et $f_1$ se rel\`event en des morphismes $\wt f_0$ et
$\wt f_1$ de $\pi_1(S)$ dans $\Z$. Ceux-ci sont en dualit\'e de Poincar\'e avec
des \'el\'ements $m_0u_0$ et $m_1u_1$ de $H_1(S;\Z)$, o\`u $m_0$, $m_1$ sont des
entiers et $u_0$, $u_1$ des classes primitives. On se donne alors des courbes
ferm\'ees simples $A_0$ et $A_1$ dont les classes d'homologie respectives sont
$u_0$ et $u_1$. Comme il existe un diff\'eomorphisme~$\phi$ de~$S$ qui envoie
$A_0$ sur~$A_1$, on peut supposer, quitte \`a composer~$f_1$ par~$\phi_*$, que
$A_0 = A_1 = A$. Le morphisme $f_i$ associe alors \`a la classe de chaque courbe
$C$ le nombre $m_i [A] \cdot [C] \mod n$. Par suite, $f_0$ et $f_1$ ont m\^eme
image dans $\Z/n\Z$ si et seulement si $m_0\Z+n\Z = m_1\Z+n\Z$, c'est-\`a-dire
si et seulement si $\pgcd(m_0,n) = \pgcd(m_1,n)$. On suppose que c'est le cas,
on note~$d$ ce plus grand commun diviseur et on pose $m_i=dm'_i$, $n=dn'$. On
choisit dans~$S$ une sous-surface compacte~$R$ contenant~$A$ et diff\'eomorphe
\`a un tore trou\'e. On note $B \subset R$ une courbe ferm\'ee simple qui, avec
$A$, forme une base de $H_1(R;\Z)$. Modulo~$n$, l'intersection de $m_i[A]$ avec
la classe d'une courbe quelconque~$C$ de~$S$ est \'egale \`a celle de $m_i [A] +
n [B]$ avec~$[C]$. Or $m_i [A] + n [B] = d (m'_i [A] + n' [B])$ et, comme
$\pgcd(m'_0,n') = \pgcd(m'_1,n') = 1$, il existe un diff\'eomorphisme~$\psi$ de
$S$, \`a support dans~$R$, qui envoie $m'_0[A] + n'[B]$ sur $m'_1[A] + n'[B]$.
Ainsi, $f_0 = f_1 \circ \psi_*$.
\end{proof}

\section{\'Etude des structures de contact non transversales}

\subsection{Structures invariantes}

On consid\`ere ici une vari\'et\'e~$V$ close, connexe et orient\'ee, munie d'une
action libre du cercle~$\S^1$ et ainsi fibr\'ee au-dessus de la surface quotient
$S = V/\S^1$. On s'int\'eresse sur~$V$ aux structures de contact invariantes par
l'action. Pour une telle structure~$\xi$, on note~$\Gamma(\xi)$ l'ensemble des
orbites $q \in S$ qui sont tangentes \`a~$\xi$ (\emph{i.e.} legendriennes). Il
est facile de voir que $\Gamma(\xi)$ est une multi-courbe lisse sur~$S$ et
R.~{\sc Lutz} montre dans~\cite{Lu} que deux structures invariantes $\xi_0$
et~$\xi_1$ sont conjugu\'ees par un diff\'eomorphisme \'equivariant de~$V$ si et
seulement s'il existe un diff\'eomorphisme de~$S$ qui envoie $\Gamma(\xi_0)$
sur~$\Gamma(\xi_1)$. Avant d'expliquer comment affranchir ces r\'esultats des
conditions d'invariance et d'\'equivariance, on \'etablit une caract\'erisation
des structures invariantes (universellement) tendues.

\begin{proposition} \label{p:lutz}
Soit $\xi$ une structure de contact orientable et invariante sur~$V$.

\alinea{\textbf{a)}}
Si $\xi$ est tendue et si une composante connexe de $S \setminus \Gamma(\xi)$ 
est un disque, $\Gamma(\xi)$ est connexe et $\Chi(V,S)$ v\'erifie
l'in\'egalit\'e
$$ \left\{ \begin{aligned}
\Chi(V,S) & > 0 &\quad& \text{si \ $S \ne \S^2$,} \\ 
\Chi(V,S) & \ge 0 &\quad& \text{si \ $S = \S^2$.}
\end{aligned} \right. $$

\alinea{\textbf{b)}}
Pour que $\xi$ soit universellement tendue, il faut et il suffit que l'une des
conditions suivantes soit remplie\up{ :}
\begin{itemize}
\item
$S \not\simeq \S^2$ et aucune composante connexe de $S \setminus \Gamma(\xi)$
n'est un disque\up{ ;}
\item
$S \simeq \S^2$, \ $\chi(V,S) < 0$ et $\Gamma(\xi)$ est vide\up{ ;}
\item
$S \simeq \S^2$, \ $\chi(V,S) \ge 0$ et $\Gamma(\xi)$ est connexe mais pas vide.
\end{itemize}
\end{proposition}

\begin{proof}
\Alinea{a)}
Soit $D$ une composante de $S \setminus \Gamma(\xi)$ qui est un disque et $E$ la
composante voisine. Si $\Gamma(\xi)$ n'est pas connexe, la surface compacte $R =
\Adh (D \cup E)$ diff\`ere de~$S$ et la fibration $\pi \from V \to S$ admet
une section $\wh R$ au-dessus de~$R$ ayant, pour un choix convenable
d'orientations, les propri\'et\'es suivantes :
\begin{itemize}
\item[1)]
la courbe $\partial \wh R$ est positivement transversale \`a~$\xi$ ;
\item[2)]
la surface~$\wh R$ a un seul point de contact n\'egatif avec~$\xi$ et ce point
est une singularit\'e d'indice~$1$ du feuilletage $\xi \wh R$.
\end{itemize}
Or l'in\'egalit\'e de Bennequin~\cite{El:martinet} interdit l'existence d'une
telle surface~$\wh R$ si $\xi$ est tendue.

On \'etablit maintenant l'in\'egalit\'e sur le nombre d'Euler. Si $\Chi(V,S) < 0$
et si $Q \subset S \setminus D$ est un disque assez petit, la fibration $\pi$
admet sur $R = S \setminus \Int Q$ une section~$\wh R$ ayant elle aussi, pour un
choix convenable d'orientations, les propri\'et\'es 1) et 2) ci-dessus qui sont
illusoires si $\xi$ est tendue. Enfin, si $\Chi(V,S) = 0$ et si $S \ne \S^2$,
le lemme~2.6 de \cite{Gi:bourbaki} (qui sert \`a \'etablir l'in\'egalit\'e de
Bennequin relative aux surfaces closes) montre encore que $\xi$ est vrill\'ee.

\smallskip

\alinea{b)}
Pour voir que l'une des conditions \'enum\'er\'ees est remplie quand $\xi$ est
universellement tendue, il suffit d'appliquer~a) et d'observer que, si $S \not
\simeq \S^2$ et si une composante connexe de $S \setminus \Gamma(\xi)$ est un
Disque, l'image inverse de $\Gamma(\xi)$ par n'importe quel rev\^etement non
trivial $\rho \from \wt S \to S$ est non connexe. Or cette multi-courbe n'est
autre que $\Gamma(\wt\xi)$ o\`u $\wt\xi$ d\'esigne le rappel de~$\xi$ sur $\wt V
= \rho^*V$. Ainsi, $\wt\xi$ est vrill\'ee et $\xi$ l'est virtuellement.

On explique maintenant pourquoi $\xi$ est universellement tendue lorsque
$\Gamma(\xi)$ satisfait l'une des conditions requises.

Si $S = \S^2$, la classification de Lutz montre que le rev\^etement universel
$(\wt V,\wt\xi)$ de $(V,\xi)$ est fait comme suit, \`a un isomorphisme pr\`es :
\begin{itemize}
\item
si $\Chi(V,S) = 0$, alors $\wt V = \S^2 \times \R = \R^3 \setminus\{0\}$ et $\wt
\xi$ est la structure usuelle, d'\'equation $dz + x\,dy - y\,dx = 0$, qui est
invariante par l'action du flot $(x,y,z) \mapsto (e^t x, e^t y, e^{2t} z)$ ;
\item
si $\Chi(V,S) = \mp1$, alors $\wt V = \S^3$ est la sph\`ere unit\'e de $\C^2$
et $\wt\xi$ est la structure usuelle, d'\'equation $\Im (\ol z\,dz + \ol w\,dw)
= 0$, qui est invariante par l'action du flot $(z,w) \mapsto (e^{i\theta} z,
e^{\pm i\theta} w)$.
\end{itemize}
Dans tous ces cas, le th\'eor\`eme de Bennequin assure que $\wt\xi$ est tendue.

Si $S \ne \S^2$, le rev\^etement universel de $S$ est $\R^2$ et il suffit de
voir
que la structure $\wt\xi$ induite par~$\xi$ sur $\wt V = \R^2 \times \S^1$ est
tendue. Comme toutes les composantes de $\Gamma(\xi)$ sont essentielles sur~$S$
(\emph{i.e.} non contractiles), celles de $\Gamma(\wt\xi)$ sont des droites
proprement plong\'ees dans~$\R^2$. On remplit $\R^2$ avec une suite exhaustive
de disques ferm\'es~$D_n$ dont les bords sont transversaux \`a $\Gamma(\wt\xi)$.
On va montrer que $\wt\xi$ est tendue en plongeant chaque domaine $(D_n \times
\S^1, \wt\xi)$ dans $(\S^2 \times \S^1, \eta)$ o\`u $\eta$ est une structure de
contact invariante ayant une courbe $\Gamma(\eta)$ connexe.

On se donne des \'equations invariantes de $\wt\xi$ et $\eta$ qu'on \'ecrit
respectivement $\beta + u\,dt = 0$ et $\lambda + v\,dt = 0$, o\`u $t \in \S^1$
et o\`u $\beta$, $u$ (resp. $\lambda$, $v$) sont une $1$-forme et une fonction
sur~$\R^2$ (resp. sur~$\S^2$). Les ensembles $\Gamma(\wt\xi)$ et $\Gamma(\eta)$
ont donc pour \'equations respectives $u=0$ et $v=0$. On choisit, pour tout $n
\ge 0$, un plongement $\phi_n \from D_n \to \S^2$ qui envoie $\Gamma(\wt\xi)
\cap D_n$ sur $\Gamma(\eta) \cap \phi_n(D_n)$ en respectant les coorientations
induites par $u$ et~$v$. Il existe ainsi une fonction $h_n \from D_n \to \op]0,
\infty\cl[$ telle que $v \circ \phi_n = h_n u$ et on pose $\beta_n = h_n \beta$.
Le lemme~\ref{l:martinet} ci-dessous garantit alors que la forme
$(\phi_n)_*\beta_n$ se prolonge \`a~$\S^2$ en une forme $\lambda_n$ v\'erifiant
l'in\'egalit\'e $v\,d\lambda_n + \lambda_n \wedge dv > 0$, laquelle assure que
l'\'equation $\lambda_n + v\,dt = 0$ d\'efinit une structure de contact $\eta_n$
invariante sur $\S^2 \times \S^1$. Or, par construction,
$$ \phi_n \times \id \from \bigl( D_n \times \S^1, \wt\xi \bigr)
   \longrightarrow \bigl( \S^2 \times \S^1, \eta_n \bigr) $$
est un plongement de contact et, comme $\Gamma(\eta_n) =\{v=0\}= \Gamma(\eta)$,
la structure $\eta_n$ est isotope \`a~$\eta$. Par cons\'equent, la structure
$\wt\xi$ est tendue et $\xi$ l'est universellement.
\end{proof}

\begin{lemme} \label{l:martinet}
Soit $S$ une surface compacte orient\'ee, $R$ une sous-surface compacte et $v
\from S \to \R$ une fonction qui admet $0$ pour valeur r\'eguli\`ere, de m\^eme
que $v \res {\partial S}$ et $v \res {\partial R}$. Si $v$ s'annule dans chaque
composante de $S \setminus R$, toute $1$-forme $\lambda$ sur $R$ qui satisfait
\`a l'in\'egalit\'e
\begin{equation*} \label{*}
  v \, d\lambda + \lambda \wedge dv > 0 \tag{$*$}
\end{equation*}
se prolonge \`a~$S$ en une $1$-forme v\'erifiant partout l'in\'egalit\'e~\eqref
{*}.
\end{lemme}

\begin{proof}
En un point de $\Gamma = \{v=0\}$, l'in\'egalit\'e~\eqref{*} dit simplement que
$\lambda$ est transversale \`a~$dv$. On prolonge donc sans peine $\lambda$ \`a
un voisinage~$U$ de $\Gamma$. D'autre part, en tout point de $S \setminus
\Gamma$, 
$$ v \, d\lambda + \lambda \wedge dv = v^2 \, d(\lambda/v) \,. $$
On observe alors que, par hypoth\`ese, chaque composante~$D$ de  $S \setminus 
(R \cup \Gamma)$ contient au moins un arc~$J$ de $\Gamma$ dans sa fronti\`ere.
Par suite, l'int\'egrale de $\lambda/v$ sur le bord de $D$ est infinie. Quitte
\`a diminuer le voisinage~$U$, on peut donc prolonger $\lambda/v$ \`a~$D$ en une
$1$-forme dont la diff\'erentielle ext\'erieure soit partout positive.
\end{proof}

\subsection{Comment revisiter Lutz}

On donne ici une description de toutes les structures de contact universellement
tendues sur une vari\'et\'e fibr\'ee en cercles au-dessus d'une surface.

\begin{definition} \label{d:cloisons}
Soit $V$ une vari\'et\'e orient\'ee fibr\'ee en cercles au-dessus d'une surface
compacte~$S$.
Une \emph{multi-courbe} sur~$S$ est ici une union disjointe d'un nombre fini de
courbes ferm\'ees simples et d'arcs proprement plong\'es dans~$S$. Par ailleurs,
une multi-courbe est \emph{essentielle} si aucune de ses composantes n'est nulle
en homotopie --~relative au bord s'il s'agit d'un arc.

On dira qu'une structure de contact $\xi$ sur~$V$ est \emph{cloisonn\'ee} par
une multi-courbe $\Gamma \subset S$ si les conditions suivantes sont remplies :
\begin{itemize}
\item
sur $V \setminus \pi^{-1}(\Gamma)$,  la structure $\xi$ est transversale aux
fibres ;
\item
la surface $\pi^{-1}(\Gamma)$ est transversale \`a $\xi$ et ses
caract\'eristiques sont des fibres.
\end{itemize}
\end{definition}

\smallskip

\begin{exemple}
Toute vari\'et\'e orient\'ee~$V$ fibr\'ee en cercles au-dessus d'une surface~$S$
(orientable) peut \^etre munie d'une action libre du cercle qui d\'efinit la
fibration. Quand $S$ est compacte, R.~{\sc Lutz} construit dans~\cite{Lu}, pour
toute multi-courbe non vide $\Gamma$ dans~$S$, une structure de contact
invariante~$\xi$ sur~$V$ telle que $\Gamma(\xi)$ soit \'egal \`a $\Gamma$. Cette
structure est alors cloisonn\'ee par~$\Gamma$. En outre, elle est orientable si
et seulement si la classe de $\Gamma$ dans $H_1(S;\Z/2\Z)$ est nulle.
\end{exemple}

\smallskip

Un th\'eor\`eme de recollement d\^u \`a V.~{\sc Colin}~\cite{Co:recollement}
assure que toute structure de contact cloisonn\'ee par une multi-courbe 
essentielle est universellement tendue\footnote
{Une autre d\'emonstration de ce fait s'ensuit de la proposition~\ref{p:lutz} et
de la partie~b) du th\'eor\`eme~\ref{t:lutz}.}.
R\'eciproquement :

\begin{theoreme} \label{t:lutz}
Soit $V$ une vari\'et\'e connexe et orient\'ee, fibr\'ee en cercles au-dessus
d'une surface close~$S$ de caract\'eristique d'Euler n\'egative ou nulle.

\alinea{\textbf{a)}}
Toute structure de contact orientable et universellement tendue est isotope \`a
une structure cloisonn\'ee par une multi-courbe essentielle.

\alinea{\textbf{b)}}
Soit $\xi_0$ et $\xi_1$ des structures de contact cloisonn\'ees par des
multi-courbes essentielles non vides, respectivement not\'es $\Gamma_{\!0}$
et~$\Gamma_{\!1}$. Les structures $\xi_0$ et~$\xi_1$ sont isotopes si et
seulement si les multi-courbes $\Gamma_{\!0}$ et~$\Gamma_{\!1}$ le sont.
\end{theoreme}

Avec les th\'eor\`emes \ref{t:ghys} et~\ref{t:thurston} --~ce dernier montrant
en particulier qu'une structure transversale aux fibres ne peut \^etre isotope
\`a une structure cloisonn\'ee par une multi-courbe non vide~--, le
th\'eor\`eme ci-dessus \'etablit une classification des structures de contact
universellement tendues sur~$V$ lorsque la caract\'eristique d'Euler de~$S$ est
n\'egative ou nulle. Pour les vari\'et\'es fibr\'ees en cercles au-dessus de la
sph\`ere, qui sont des espaces lenticulaires, le th\'eor\`eme~1.1 de~\cite
{Gi:bifurcations} donne une classification de toutes les structures de contact
tendues.

\medskip

D'autre part, la partie~a) du th\'eor\`eme~\ref{t:lutz} prouve que toute
structure de contact universellement tendue et d'enroulement positif ou nul est
isotope \`a une structure invariante (par une quelconque action libre du cercle
qui d\'efinit la fibration). La partie~b) classifie donc en fait les structures
de contact $\S^1$-invariantes. Sa d\'emonstration s'adapte alors, sans surprise,
aux structures de contact $\R$-invariantes sur le produit d'une surface par
$\R$. On obtient ainsi, compte tenu de l'abondance des surfaces convexes (cf.
section 2.D), une classification \guil{g\'en\'erique} des structures de contact
tendues au voisinage des surfaces :

\begin{theoreme} \label{t:convexe}
Soit $(M,\xi)$ une vari\'et\'e de contact de dimension~$3$, $F \subset (M,\xi)$
une surface convexe close, $U = F \times \R$ un voisinage homog\`ene de $F$ et
$\Gamma$ la multi-courbe qui scinde $\xi F$.

\alinea{\textbf{a)}}
La restriction de $\xi$ \`a $U$ est tendue si et seulement si l'une des
conditions suivantes est remplie\up{ :}
\begin{itemize}
\item
$F \not\simeq \S^2$ et aucune composante de $F \setminus \Gamma$ n'est un
disque\up{ ;}
\item
$F \simeq \S^2$ et $\Gamma$ est connexe mais pas vide.
\end{itemize}

\alinea{\textbf{b)}}
On suppose que $\xi$ est tendue. Une surface convexe $F' \subset (M,\xi)$
poss\`ede un voisinage homog\`ene isomorphe \`a~$(U,\xi)$ si et seulement s'il
existe un diff\'eomorphisme de $F$ dans $F'$ qui envoie $\Gamma$ sur une
multi-courbe qui scinde le feuilletage caract\'eristique $\xi F'$.
\end{theoreme}

La partie~a) de ce th\'eor\`eme est un corollaire imm\'ediat de la proposition
\ref{p:lutz}. En effet, la restriction de $\xi$ \`a $U = F \times \R$ est tendue
si et seulement si la structure de contact $\ol\xi$ induite par~$\xi$ sur $F
\times \R/n\Z$ est tendue pour tout entier $n>0$. En outre, la multi-courbe
$\Gamma_{\!U}$ (cf.~d\'efinition~\ref{d:ctc}) qui scinde le feuilletage $\xi F$
n'est autre que $\Gamma(\ol\xi)$. La partie~b) sera d\'emontr\'ee dans la
section~D.

\subsection{Existence d'un cloisonnement}

On d\'emontre ici la partie~a) du th\'eor\`eme~\ref{t:lutz}. La surface~$S$ est
donc de caract\'eristique d'Euler n\'egative ou nulle et les structures de
contact qu'on consid\`ere sont orientables.

\begin{lemme} \label{l:roussarie'}
Soit $\xi$ une structure de contact tendue sur~$V$ et $R \subset S$ une surface
compacte, connexe et \`a bord non vide. Si l'enroulement $\w(\xi)$ de~$\xi$ est
positif ou nul, $\xi$ est isotope \`a une structure qui, au-dessus de~$R$, est
cloisonn\'ee par un syst\`eme d'arcs.
\end{lemme}

\begin{proof}
On adapte la d\'emonstration du lemme~\ref{l:roussarie}. On regarde $R$ comme un
voisinage r\'egulier d'un bouquet de cercles~$K$ de sommet~$q$. Quitte \`a faire
une premi\`ere isotopie, on peut trouver des coordonn\'ees $(x,y,t) \in \D^2
\times \S^1$, au-dessus d'un voisinage compact $Q$ de~$q$, dans lesquelles~$\pi$
est la projection sur le disque, $\xi$ a pour \'equation $dy+x\,dt = 0$ et $q =
(0,0)$. Sur $N = \pi^{-1}(Q)$, la structure~$\xi$ est cloisonn\'ee par l'arc
$J = \{x=0\}$.

On modifie maintenant~$K$, par des mouvements de Whitehead \`a support dans~$Q$,
en un graphe~$K'$ constitu\'e de cercles lisses $K'_i$, $1 \le i \le k$, et d'un
arbre $K'_0$ inclus dans~$Q$. On met en outre tous les sommets de $K'$ sur l'arc
$J \subset Q$ et on rend les ar\^etes transversales d'une part \`a~$J$, d'autre
part au vecteur $\partial_y$ en tout point d'intersection avec~$J$. On d\'eforme
ensuite~$\xi$, par une petite isotopie relative \`a~$N$, pour rendre convexe la
surface $\pi^{-1} (K' \setminus \Int Q)$, qui est une union d'anneaux. Chaque
tore $F_i = \pi^{-1}(K'_i)$, $1 \le i \le k$, est ainsi convexe et la
multi-courbe qui
scinde son feuilletage caract\'eristique $\xi F_i$ est verticale. Elle ne peut
en effet intersecter les fibres au-dessus de $K'_i \cap J$ qui sont des feuilles
ferm\'ees de $\xi F_i$. On consid\`ere alors sur~$F_i$ un feuilletage singulier
$\sigma_i$ ayant les propri\'et\'es suivantes :
\begin{itemize}
\item
$\sigma_i$ est scind\'e par la m\^eme multi-courbe que $\xi F_i$ ;
\item
$\sigma_i$ co\"{\i}ncide avec $\xi F_i$ dans $F_i \cap \pi^{-1}(Q)$ ;
\item
$\sigma_i$ est non singulier et est transversal aux fibres en dehors de ses
feuilles ferm\'ees, lesquelles sont des fibres.
\end{itemize}
Le lemme~\ref{l:ctc} fournit, pour $1 \le i \le k$, un plongement $\phi_i$ de
$F_i$ dans~$V$ --~isotope \`a l'inclusion~-- dont l'image a pour feuilletage
caract\'eristique $(\phi_i)_*\sigma_i$ et a m\^eme intersection que $F_i$ avec
le compact
$$ P_i = \pi^{-1}(J \cup K'') \cup \phi_1(F_1) \cup \dots \cup
   \phi_{i-1}(F_{i-1}) \cup F_{i+1} \cup \dots \cup F_k \,. $$
Si $\phi$ est un diff\'eomorphisme de $V$ isotope \`a l'identit\'e qui prolonge
les divers plongements~$\phi_i$, la structure de contact $\xi' = \phi^*\xi$
trace sur chaque tore~$F_i$ le feuilletage~$\sigma_i$.

Pour compl\'eter la preuve, on param\`etre un voisinage tubulaire de $K'_i$ par
$K'_i \times [-1,1]$, o\`u $K'_i = K'_i \times \{0\}$. Les tores $F_{i,s} =
\pi^{-1} (K'_i \times \{s\})$ ont, pour $s$ petit, un feuilletage $\xi' F_{i,s}$
conjugu\'e \`a $\xi' F_i$ car ce feuilletage est topologiquement stable. On peut
donc redresser $\xi'$, par une petite isotopie stationnaire sur $\bigcup F_i
\cup \pi^{-1}(Q)$, en une structure $\xi''$ pour laquelle $\xi''F_{i,s}$ est
transversal aux fibres
en dehors de ses feuilles ferm\'ees, lesquelles sont des fibres. Les projections 
sur~$S$ des feuilles ferm\'ees de tous les feuilletages $\xi'' F_{i,s}$, avec
$s$ petit, forment alors un syst\`eme d'arcs qui, avec~$J$, cloisonne $\xi''$
sur un voisinage de $K'$. Comme $R$ se r\'etracte par isotopie sur un voisinage
arbitrairement petit de~$K'$, le lemme est d\'emontr\'e.
\end{proof}

\smallskip

D\'esormais, $\xi$ d\'esigne une structure de contact universellement tendue sur
$V$. On note $A$ un anneau non s\'eparant dans~$S$ et on pose $R = S \setminus
\Int A$. Compte tenu du lemme ci-dessus, on suppose que $\xi$ est cloisonn\'ee,
au-dessus de $R$, par un syst\`eme d'arcs $\Gamma_{\!R} \subset R$. On
param\`etre
$A$ par $\S^1 \times [0,1]$ et $W = \pi^{-1}(A)$ par $\T^2 \times [0,1]$ de
telle sorte que la fibration $\pi \res W$ soit la projection. Pour tout $a \in
[0,1]$, on pose encore $T_a = \T^2 \times \{a\}$. L'argument utilis\'e au lemme
\ref{l:univ} montre ici que la restriction de $\xi$ \`a $W$ est universellement
tendue. D'apr\`es les propositions 3.15, 3.22 et 3.29 de \cite{Gi:bifurcations},
la structure $\xi \res W$ est alors isotope, relativement au bord de~$W$, \`a
une structure de contact~$\eta$ pour laquelle il existe dans $A \simeq \S^1
\times [0,1]$ une multi-courbe $\Gamma_{\!A}$ ayant les propri\'et\'es
suivantes : 
\begin{itemize}
\item
si $T_a \cap \pi^{-1}(\Gamma_{\!A}) \ne \varnothing$, cette intersection est
l'union des feuilles ferm\'ees et des singularit\'es de $\eta T_a$  (ces
singularit\'es formant donc des courbes) ;
\item
si $T_a \cap \pi^{-1}(\Gamma_{\!A}) = \varnothing$, le feuilletage $\eta T_a$
est une suspension dont aucune feuille ferm\'ee n'est isotope \`a la fibre.
\end{itemize}
\`A partir de l\`a, on d\'eforme facilement~$\eta$, par une isotopie relative
\`a $\partial W \cup \pi^{-1}(\Gamma_{\!A})$, en une structure de contact~$
\eta'$ cloisonn\'ee par $\Gamma_{\!A}$. En recollant $\eta'$ avec la restriction
de~$\xi$ \`a $\pi^{-1}(R)$, on obtient une structure de contact~$\xi'$
cloisonn\'ee par $\Gamma = \Gamma_{\!R} \cup \Gamma_{\!A}$. Il reste \`a montrer
que les courbes de~$\Gamma$ sont toutes essentielles, fait qui r\'esulte de la
proposition~\ref{p:lutz} et du lemme suivant :

\begin{lemme} \label{l:lutz}
Deux structures de contact cloisonn\'ees par une m\^eme multi-courbe non vide
sont isotopes.
\end{lemme}

\begin{proof}
Soit $\xi_0$ et $\xi_1$ les structures de contact, $\Gamma$ l multi-courbe qui
les cloisonne et $u \from S \to \R$ une fonction dont le niveau $\{u = 0\}$ est
r\'egulier et \'egal \`a $\Gamma$. D'apr\`es \cite{Lu} (voir aussi
\cite{Gi:convexite}), il existe sur $S$ une $1$-forme~$\beta$ pour laquelle la
$2$-forme $u\, d\beta + \beta \wedge du$ est une forme d'aire sur~$S$.

On munit maintenant  $V$ d'une action libre du cercle d\'efinissant la fibration
et d'une forme de connexion~$\tau$. Comme chacune des structures $\xi_i$, $i \in
\{0,1\}$, est cloisonn\'ee par~$\Gamma$, elle admet une \'equation du type
$\beta_i + \pi^*u\, \tau = 0$ o\`u $\beta_i$ est une $1$-forme sur $V$ nulle sur
les vecteurs tangents aux fibres. Un calcul direct montre alors que, si $s$ est
un r\'eel positif pris assez grand, les \'equations de Pfaff
$$ (1-t)\beta_i + ts\pi^*\beta + \pi^*u\, \tau = 0, \qquad i \in \{0,1\}, $$
d\'efinissent des structures de contact pour tout $t \in [0,1]$. Le th\'eor\`eme
de Gray assure d\`es lors que $\xi_0$ et $\xi_1$ sont isotopes.
\end{proof}

\subsection{In\'egalit\'e de Bennequin semi-locale}

On d\'emontre ici la partie~b) du th\'eor\`eme \ref{t:convexe}. S'il existe un
diff\'eomorphisme de $F$ dans $F'$ envoyant $\Gamma$ sur une multi-courbe qui
scinde $\xi F'$, les r\'esultats de~\cite{Lu} assurent que tout voisinage
homog\`ene de $F'$ est isomorphe \`a $(U,\xi)$. Pour \'etablir la r\'eciproque,
on utilise la notion d'intersection g\'eom\'etrique.

\begin{definition} \label{d:thurston}
Sur une surface close, on consid\`ere une courbe ferm\'ee simple~$C$ et une
multi-courbe~$\Gamma$. L'\emph{intersection g\'eom\'etrique} $\q(\Gamma,C)$ est
le nombre minimal de points d'intersection entre $\Gamma$ et une courbe
quelconque isotope \`a~$C$.
\end{definition}

Un des int\'er\^ets de cette notion r\'eside dans la proposition suivante, qui
est \`a la base des travaux de W.~{\sc Thurston} sur les surfaces~\cite
{Th:surfaces} (voir aussi \cite[expos\'e~4 p.~59]{FLP}) :

\begin{proposition} \label{p:flp}
Soit $\Gamma_{\!0}$ et $\Gamma_{\!1}$ deux multi-courbes essentielles sur une
surface close. Si $\q(\Gamma_{\!0},C) = \q(\Gamma_{\!1},C)$ pour toute courbe
ferm\'ee simple~$C$, alors $\Gamma_{\!0}$ et $\Gamma_{\!1}$ sont isotopes.
\end{proposition}

D'apr\`es un th\'eor\`eme de J.~{\sc Stallings}, tout diff\'eomorphisme de $F
\times \R$ dans $F' \times \R$ est isotope \`a un diff\'eomorphisme produit. La
partie~b) du th\'eor\`eme~\ref{t:convexe} d\'ecoule alors directement de la
proposition ci-dessus et de l'\emph{in\'egalit\'e de Bennequin semi-locale} que
voici :

\begin{proposition} \label{p:bennequin}
Soit $\xi$ une structure de contact $\R$-invariante et tendue sur le produit $U
= F \times \R$. Soit $C$ une courbe ferm\'ee simple sur $F = F \times \{0\}$ et
$\Gamma$ une multi-courbe qui scinde $\xi F$. Pour toute isotopie $\phi_t$ de
$U$ qui am\`ene~$C$ sur une courbe legendrienne~$\phi_1(C)$, le nombre de tours
que fait $\xi$ par rapport au plan tangent \`a $\phi_1(F)$ le long de
$\phi_1(C)$ v\'erifie l'in\'egalit\'e
$$ \deg \bigl( \xi, \phi_1(F); \phi_1(C) \bigr)
   \le - \tfrac12 \q(\Gamma,C) \,. $$
De plus, il existe une isotopie~$\phi_t$ qui r\'ealise l'\'egalit\'e.
\end{proposition}

\begin{proof}
Quitte \`a d\'eplacer $C$ sur~$F$ par une isotopie, on suppose que $C$ rencontre
$\Gamma$ en $\q(\Gamma,C)$~points et transversalement. On prend en outre $\xi$
orientable --~ce qui revient \`a passer \'eventuellement \`a un rev\^etement
double~-- et, comme la classe de $\Gamma$ dans $H_1(F;\Z/2\Z)$ est alors nulle,
on pose $2n = \q(\Gamma,C)$.

On observe tout d'abord que l'\'egalit\'e est atteinte. En effet, on construit
facilement un feuilletage $\sigma$ de~$F$ scind\'e par~$\Gamma$ et dans lequel
$C$ est une union de singularit\'es et de feuilles. De plus, le lemme~\ref
{l:ctc} fournit une isotopie $\phi_t$ de plongements de $F$ dans~$U$ qui
am\`ene $F$ sur une surface $\phi_1(F)$ ayant pour feuilletage caract\'eristique
$(\phi_1)_*\sigma$. La courbe $\phi_1(C)$ est alors une courbe legendrienne le
long de laquelle $\xi$ fait $-\frac12 \Card (\Gamma \cap C)$ tours par rapport
au plan tangent \`a $\phi_1(F)$.

\smallskip

Soit $\rho \from \wt F \to F$ le rev\^etement associ\'e \`a~$C$ et $\wt C$ un
rel\`evement compact de~$C$ dans $\wt F$. On note~$\wt\xi$ le rappel de $\xi$
sur $\wt U = \wt F \times \R$ et $\wt\phi_t$, $t \in [0,1]$, le rel\`evement de
l'isotopie $\phi_t$ \`a $\wt U$. Le nombre de tours que fait $\xi$ par rapport
\`a $\phi_1(F)$ le long de $\phi_1(C)$ est clairement \'egal au nombre de tours
que fait $\wt\xi$ par rapport \`a $\wt\phi_1 (\wt F)$ le long de $\wt\phi_1 (\wt
C)$. En outre, les courbes $\wt\phi_t (\wt C)$ restent dans un compact de $\wt
U$. On param\`etre alors $\wt F$ par $\S^1 \times \R$ de telle sorte que $\wt C$
soit la courbe $\S^1 \times \{0\}$ et on se donne un r\'eel $a>0$ assez grand
pour que toutes les courbes $\wt\phi_t(\wt C)$ soient contenues dans le domaine
$\wt U_a = \wt F_a \times \R$, o\`u $\wt F_a = \S^1 \times [-a,a]$. On \'etablit
ci-dessous l'in\'egalit\'e voulue en plusieurs \'etapes. On plonge d'abord $(\wt
U_a, \wt\xi)$ dans un mod\`ele abstrait, puis on r\'ealise ce mod\`ele dans la
sph\`ere $\S^3$ munie de sa structure de contact ordinaire et on conclut \`a
l'aide de l'in\'egalit\'e de Bennequin classique.

\smallskip

\begin{assertion}
Aucune courbe de $\wt\Gamma = \rho^{-1}(\Gamma)$ ne coupe $\wt C$ en plus d'un
point.
\end{assertion}

\begin{proof}[Preuve]
On suppose qu'une courbe de $\wt\Gamma$ coupe $\wt C$ en deux points et on note
$J$ un arc de cette courbe joignant deux points d'intersection cons\'ecutifs. La
composante connexe born\'ee de $\wt F \setminus (\wt C \cup J)$ est un disque
$D$ et $D \setminus \rho^{-1}(C)$ a au moins une composante connexe~$D_0$ dont
le bord est l'union de deux arcs, l'un contenu dans~$J$ et l'autre dans
$\rho^{-1}(C)$. La restriction de $\rho$ \`a~$D_0$ est alors injective et, en
d\'epla\c cant~$C$ par isotopie le long de $\rho(D_0)$, on \'elimine deux points
d'intersection avec~$\Gamma$. Ceci est absurde puisque $\Card (\Gamma \cap C) =
\q(\Gamma,C)$.
\end{proof}

\smallskip

L'assertion ci-dessus assure que les composantes de $\wt\Gamma$ qui vont d'un
bout \`a l'autre de $\wt F$ sont exactement celles qui rencontrent $\wt C$ et
leur nombre est donc \'egal \`a $2n = \q(\Gamma,C)$. Selon que $n$ est nul ou
non, il existe alors un plongement incompressible $\psi_{a,n}$ de $\wt F_a$ dans
$\T^2$ ou dans $\S^1 \times \R$ tel que
$$ \psi_{a,n} \bigl( \wt\Gamma \cap \wt F_a \bigr)
 = \psi_{a,n} \bigl(\wt F_a \bigr) \cap \Gamma_n, $$
o\`u
\begin{align*}
\Gamma_0 & = \bigl\{ (x,y) \in \R^2/\!\Z^2 \mid y=\pm1/4 \bigr\} \quad
  \text{et} \\
\Gamma_n & = \bigl\{ (x,y) \in \R/\Z \times \R \mid nx=0 \bigr\},
  \quad \text{pour $n>0$.}
\end{align*}
Par suite, si $n=0$ (resp. si $n>0$), le lemme~\ref{l:martinet} permet, comme
dans la d\'emonstration de la proposition~\ref{p:lutz}, de plonger
incompressiblement $(\wt U_a,\wt\xi)$ dans $(\T^2 \times \R, \xi_0)$ \ (resp.
dans $(\S^1 \times \R^2, \xi_n)$) o\`u $\xi_n$ est n'importe quelle structure de
contact $\R$-invariante qui imprime sur $\T^2 \times \{0\}$ (resp. sur $\S^1
\times \R \times \{0\}$) un feuilletage caract\'eristique scind\'e par
$\Gamma_n$.

\smallskip

Soit maintenant $\zeta$ la structure de contact ordinaire sur la sph\`ere $\S^3
\subset \C^2$ et soit $L_0$ la courbe legendrienne (non nou\'ee) $\S^3 \cap
\R^2$ --~dont l'invariant de Thurston-Bennequin $\tb(L_0)$ vaut~$-1$. 

Le th\'eor\`eme de Darboux permet de param\'etrer un voisinage~$W$ de~$L_0$ par
$\S^1 \times \R^2$ de telle sorte que $L_0$ soit la courbe $\S^1 \times \{0\}$
et que $\zeta$ ait pour \'equation $dz + p\, d\theta = 0$, $(\theta,p,z) \in
\S^1 \times \R^2$. Le tore
$$ T = \bigl\{ (\theta,p,z) \in \S^1 \times \R^2
  \mid \abs{p}^2 + \abs{z}^2 = 1 \bigr\} $$
est convexe car son feuilletage caract\'eristique $\zeta T$ est scind\'e par les
deux cercles $\{p=\pm1\}$. Par suite, $T$ poss\`ede un voisinage homog\`ene
isomorphe \`a $(\T^2 \times \R, \xi_0)$. Si $\q(\Gamma,C) = 0$, on peut donc
plonger $(\wt U_a, \wt\xi)$ dans $(\S^3, \zeta)$ en envoyant $\wt\phi_1(\wt C)$
sur une courbe legendrienne non nou\'ee $L$ dont l'invariant de
Thurston-Bennequin vaut
$$ \tb(L) = \deg \bigl( \wt\xi, \phi_1(\wt F); \wt\phi_1(\wt C) \bigr) - 1 \,.$$
L'in\'egalit\'e de Bennequin assure alors que le degr\'e est n\'egatif ou nul.

Pour finir, on consid\`ere sur $\S^1 \times \R^2$ la structure de contact
$\xi_n$ d'\'equation
$$ \cos (2n\pi x)\, dy - \sin(2n\pi x)\, dt = 0 \,. $$
Un calcul direct montre que le plongement de $\S^1 \times \R^2$ dans $\S^3$
donn\'e par
$$ (x,y,t) \longmapsto \bigl( \theta = 2\pi x, \;
  z = \cos(2n\pi x)\, y - \sin(2n\pi x)\, t, \;
  p/n = \sin(2n\pi x)\, y + \cos(n\theta)\, t \bigr) $$
envoie $\xi_n$ sur~$\zeta$. Si $\q(\Gamma,C) = 2n > 0$, on peut donc plonger,
par composition, $(\wt U_a, \wt\xi)$ dans $(\S^3, \zeta)$ en envoyant
$\wt\phi_1(\wt C)$ sur une courbe legendrienne non nou\'ee~$L$ dont l'invariant
de Thurston-Bennequin vaut
$$ \tb(L)
 = \deg \bigl( \wt\xi, \wt\phi_1(\wt F); \wt\phi_1(\wt C) \bigr) + n -1 \,. $$
L'in\'egalit\'e de Bennequin assure alors que le degr\'e vaut au plus~$-n$.
\end{proof}

\subsection{Unicit\'e du cloisonnement.}

On d\'emontre ici la partie~b) du th\'eor\`eme \ref{t:lutz}. Compte tenu du
lemme~\ref{l:lutz}, il suffit de prouver que, si les structures de contact
$\xi_0$ et~$\xi_1$ sont isotopes, les multi-courbes qui les cloisonnent le sont
aussi. L'argument est une variante de celui qui conduit \`a l'in\'egalit\'e de
Bennequin semi-locale. Il repose sur une interpr\'etation appropri\'ee de
l'intersection g\'eom\'etrique. Par commodit\'e, on appelle dans la suite \emph
{indice} d'un tore convexe dans une vari\'et\'e de contact tendue le nombre de
composantes connexes de toute multi-courbe qui scinde son feuilletage
caract\'eristique.

\begin{lemme} \label{l:indice}
Soit $\xi$ une structure de contact sur~$V$ cloisonn\'ee par une multi-courbe
essentielle~$\Gamma$ et soit~$C$ une courbe ferm\'ee simple de~$S$ dont
l'intersection g\'eom\'etrique avec $\Gamma$ n'est pas nulle. L'indice minimal
des tores convexes isotopes \`a $\pi^{-1}(C)$ est \'egal \`a $\q(\Gamma,C)$ et
leur feuilletage caract\'eristique est scind\'e par des courbes isotopes aux
fibres.
\end{lemme}

\begin{proof}
On suppose que $C$ intersecte $\Gamma$ en $\q(\Gamma,C)$ points et on note $F_0$
le tore $\pi^{-1}(C)$. Au-dessus de $C \setminus \Gamma$ (resp. de $\Gamma \cap
C$), le feuilletage $\xi F_0$ est transversal (resp.  tangent) aux fibres. En
particulier, chaque fibre $\pi^{-1}(q)$, $q \in \Gamma \cap C$, est une feuille
ferm\'ee ou une courbe de singularit\'es de $\xi F_0$. Comme $\q(\Gamma,C) \ne
0$, le tore~$F_0$ est convexe et son indice est \'egal \`a $\q(\Gamma,C)$.

Dans la suite, on d\'esigne par $\rho \from \wt S \to S$ le rev\^etement
associ\'e \`a la courbe~$C$, par $\wt\pi \from \wt V = \rho^*V \to \wt S$ la
fibration induite et par $\wt\xi$ le rappel de $\xi$ sur $\wt V$.
On pose $\wt\Gamma = \rho^{-1}(\Gamma)$ et $\wt F_0 = \wt\pi^{-1} (\wt C)$
o\`u $\wt C$ est un rel\`evement compact de~$C$ dans $\wt S$. En outre, on
param\`etre $\wt S$ par $\R \times \S^1$ de telle sorte que $\wt C$ soit $\{0\}
\times \S^1$ et $\wt V$ par $\R \times \T^2$ de telle sorte que $\wt\pi$ soit la
projection.

Soit maintenant $F$ un tore convexe isotope \`a $\pi^{-1}(C)$ et $\wt F$ le
rel\`evement compact de~$F$ dans $\wt V$ obtenu en relevant depuis $\wt F_0$ une
isotopie entre $F_0$ et~$F$. Soit encore $a$ et $\eps$ des r\'eels positifs
satisfaisant aux conditions suivantes :
\begin{itemize}
\item
les cercles $\{s\} \times \S^1$, $s \in [-\eps,\eps]$, sont tous transversaux
\`a $\wt\Gamma$ --~et coupent donc $\wt\Gamma$ en $\q(\Gamma,C)$~points ;
\item
les cercles $\{\pm a\} \times \S^1$ sont transversaux \`a $\wt\Gamma$ et le
domaine $[-a,a] \times \T^2 \subset \wt V$ contient $\wt F$.
\end{itemize}
Comme dans la d\'emonstration de la proposition~\ref{p:bennequin}, les courbes
de $\wt\Gamma$ qui vont d'un bord \`a l'autre de l'anneau $[-a,a] \times \S^1$
sont celles qui intersectent~$\wt C$ et sont en nombre~$\q(\Gamma,C)$. Par
suite, il existe un plongement
$$ \psi \from [-a,a] \times \S^1 \longrightarrow [-\eps,\eps] \times \S^1 $$
qui est l'identit\'e sur $\{0\} \times \S^1$ et v\'erifie
$$ \psi \Bigl( \wt\Gamma \cap \bigl( [-a,a] \times \S^1 \bigr) \Bigr)
 = \wt\Gamma \cap \psi \bigl( [-a,a] \times \S^1 \bigr) \,. $$
Le lemme~\ref{l:martinet} fournit alors, comme dans la proposition~\ref{p:lutz},
un plongement de contact
$$ \phi \from \bigl( [-a,a] \times \T^2, \wt\xi \bigr) \longrightarrow
   \bigl( [-\eps,\eps] \times \T^2, \wt\xi \bigr) $$
qui induit l'identit\'e sur $\wt F_0 = \{0\} \times \T^2$. En outre,
l'in\'egalit\'e de Bennequin semi-locale (proposition~\ref{p:bennequin}) montre
que le feuilletage caract\'eristique d'un tore convexe isotope \`a $\{0\} \times
\T^2$ dans $([-\eps,\eps] \times \T^2, \wt\xi)$ est scind\'e par au moins
$\q(\Gamma,C)$~courbes qui sont toutes isotopes aux fibres. En particulier,
l'indice de~$F$ --~qui est \'egal \`a celui de $\wt F$ donc \`a celui de
$\phi(\wt F)$~-- vaut au moins $\q(\Gamma,C)$.
\end{proof}

\smallskip

On compl\`ete \`a pr\'esent la d\'emonstration du th\'eor\`eme~\ref{t:lutz}-b.
Soit $C$ une courbe ferm\'ee simple sur $S$. Compte tenu de la proposition~\ref
{p:flp} et du lemme~\ref{l:indice}, il suffit de montrer que, lorsque
$\q(\Gamma_{\!0},C)$ est non nul, $\q(\Gamma_{\!1},C)$ l'est aussi. On suppose
donc dans la suite que $\q(\Gamma_{\!0},C) \ne 0$ tandis que $\q(\Gamma_{\!1},C)
= 0$.

Soit $C_i$, $i \in \{0,1\}$, une courbe isotope \`a $C$ qui intersecte
$\Gamma_{\!i}$ en $\q(\Gamma_{\!i},C)$~points et soit $F_i$ le tore
$\pi^{-1}(C_i)$. Comme $\q(\Gamma_{\!0},C) \ne 0$, le tore $F_0$ est convexe et
son indice est \'egal \`a $\q(\Gamma_{\!0},C)$. D'autre part, comme
$\q(\Gamma_{\!1},C) = 0$, le feuilletage $\xi_1 F_1$ est transversal aux fibres.
Si on perturbe $F_1$ en un tore convexe~$F$ par une isotopie assez petite, le
feuilletage caract\'eristique $\xi_1 F$ est scind\'e par des courbes qui ne sont
pas isotopes aux fibres. Il r\'esulte alors du lemme~\ref{l:indice} que les
structures de contact $\xi_0$ et $\xi_1$ ne sont pas isotopes.
\qed

\subsection{Structures virtuellement vrill\'ees}

On termine cet expos\'e par un r\'esultat de finitude pour les structures de
contact virtuellement vrill\'ees.

\begin{theoreme} \label{t:bennequin}
Soit $V$ une vari\'et\'e connexe et orient\'ee, fibr\'ee en cercles au-dessus
d'une surface close~$S$. Les structures de contact orientables et virtuellement
vrill\'ees sur $V$ forment un nombre fini de classes d'isotopie born\'e par
$$ \left\{ \begin{aligned}
    \sup \bigl\{ 0, - \Chi(S) - \Chi(V,S) - 1 \bigr\} && \quad
& \text{si $\Chi(V,S) \le 0$,} \\
1 + \sup \bigl\{ 0, - \Chi(S) - \Chi(V,S) - 1 \bigr\} && \quad
& \text{si $\Chi(V,S) > 0$.}
\end{aligned} \right. $$
\end{theoreme}

En fait, avec les formes normales d\'egag\'ees dans \cite{Gi:bifurcations}, on
obtient une description pr\'ecise de tous les exemples potentiels de (classes
d'isotopie de) structures de contact virtuellement vrill\'ees sur~$V$. Il est
par
ailleurs probable que toutes les structures de contact ainsi d\'ecrites sont
effectivement tendues --~m\^eme holomorphiquement remplissables~-- et que les
techniques de chirurgie d\'evelopp\'ees par R.~{\sc Gompf} dans \cite{Go:stein}
permettraient de le prouver. Du reste, la proposition~\ref{p:noeuds} montre, par
des astuces de rev\^etements, l'existence de structures de contact virtuellement
vrill\'ees sur les vari\'et\'es fibr\'ees en cercles au-dessus du tore et dont
le nombre d'Euler est inf\'erieur ou \'egal \`a~$-2$.

\begin{lemme} \label{l:virt}
Si $V$ porte une structure de contact virtuellement vrill\'ee et d'enroulement
positif ou nul, celle-ci est isotope \`a une structure cloisonn\'ee par une
courbe connexe contractile et $\Chi(V,S)$ est strictement positif.
\end{lemme}

\begin{proof}
Soit $A \subset S$ un anneau, $R$ la surface $S \setminus \Int A$ et $W$ le tore
\'epais $\pi^{-1}(A)$. D'apr\`es le lemme \ref{l:roussarie'}, toute structure de
contact tendue d'enroulement positif ou nul est isotope \`a une structure~$\xi$
qui, au-dessus de $R$, est cloisonn\'ee par un syst\`eme d'arcs~$\Gamma_{\!R}$.
La restriction de $\xi$ \`a $W \simeq \T^2 \times [0,1]$ est ainsi une structure
de contact tendue qui trace sur chaque composante de $\partial W$ un feuilletage
ayant des feuilles ferm\'ees ou des cercles de singularit\'es parall\`eles aux
fibres. Vu ce comportement au bord, le th\'eor\`eme 1.5 de~\cite
{Gi:bifurcations} assure que $\xi \res W$ est universellement tendue et est
isotope, relativement \`a $\partial W$, \`a une structure de contact~$\eta$
cloisonn\'ee par un syst\`eme d'arcs~$\Gamma_{\!A}$. Ainsi $\xi$ est isotope \`a
une structure de contact cloisonn\'ee par la multi-courbe~$\Gamma = \Gamma_{\!R}
\cup \Gamma_{\!A}$.

Pour conclure, on se donne sur~$V$ une action libre du cercle dont les orbites
sont les fibres de~$\pi$. D'apr\`es le lemme~\ref{l:lutz}, $\xi$ est isotope \`a
\guil{la} structure invariante cloisonn\'ee par~$\Gamma$. Comme $\xi$ est tendue
mais pas universellement, la proposition~\ref{p:lutz} montre que $\Gamma$ est
une courbe connexe contractile et que $\Chi(V,S) > 0$.
\end{proof}

\smallskip

\begin{proof}[D\'emonstration du th\'eor\`eme~\ref{t:bennequin}]
D'apr\`es les lemmes \ref{l:virt}, \ref{l:lutz} et la proposition~\ref{p:lutz},
les structures de contact virtuellement vrill\'ees et d'enroulement positif ou
nul forment au plus une classe d'isotopie, et aucune si $\Chi(V,S) \le 0$. On
\'etudie donc d\'esormais les structures de contact virtuellement vrill\'ees et
d'enroulement strictement n\'egatif.

Soit $D \subset S$ un disque, $R$ la surface $S \setminus \Int D$ et $W$ le tore
plein $\pi^{-1}(D)$. Comme d'habitude, on param\`etre $W$ par $\D^2 \times \S^1$
de telle sorte que $\pi \res W$ soit la projection sur~$\D^2$. D'apr\`es le
lemme~\ref{l:roussarie}, toute structure de contact d'enroulement $-n$, $n>0$,
est isotope \`a une structure~$\xi$ pour laquelle les fibres au-dessus de~$R$
sont legendriennes et d'enroulement~$-n$. Dans ces conditions, $\partial W$ est
un tore convexe d'indice~$2$ (lemme~\ref{l:deux}) et les singularit\'es de son
feuilletage caract\'eristique $\xi\, \partial W$ forment des courbes de classe
$(n, n\Chi(V,S) + \Chi(S) - 1)$ dans $H_1(\partial W;\Z) \cong \Z^2$ (corollaire
\ref{c:whitney}).

\begin{assertion}
Ou bien $n \Chi(V,S) = -\Chi(S)$, ou bien $n=1$ et $\Chi(V,S) < -\Chi(S)$.
\end{assertion}

\begin{proof}[Preuve]
Pour tout $a \in \op]0,1]$, on note~$T_a$ le tore $a\S^1 \times \S^1$. D'apr\`es
la proposition~3.22 de~\cite{Gi:bifurcations}, la restriction de $\xi$ \`a~$W$
est isotope, relativement au bord, \`a une structure de contact $\eta$
transversale \`a $\{0\} \times \S^1$ et dont les feuilletages caract\'eristiques
$\eta T_a$ ont les propri\'et\'es suivantes :
\begin{itemize}
\item
$\eta T_a$ est une suspension sauf pour un nombre fini de valeurs $a_1, \dots,
a_k \in \op]0,1]$ ;
\item
$\eta T_{a_i}$, $1\le i \le k$, n'a aucune orbite ferm\'ee et ses singularit\'es
forment deux cercles.
\end{itemize}
Chaque feuilletage $\eta T_a$ d\'etermine alors une droite $\delta_a$ dans $\R^2
\cong H_1(T_a,\R)$ qui, pour $a \notin \{a_i\}$, porte les cycles asymptotiques
et, pour $a \in \{a_i\}$, contient la classe des cercles de singularit\'es.
Cette droite varie contin\^ument avec~$a$ et converge vers $\delta_0 = \R \times
\{0\}$ quand $a$ tend vers~$0$. Les droites $\delta_a$, $a\in[0,1]$, d\'ecrivent
donc un connexe~$\Delta$ de $\RP^1$, connexe qui ne contient pas la droite
$\{0\} \times \R$ car l'enroulement de $\eta$ est strictement n\'egatif. Ainsi,
$\Delta$ est l'intervalle $[\delta_1,\delta_0]$ pour l'orientation naturelle de
$\RP^1$ et $n \Chi(V,S) + \Chi(S) \le 0$ puisque $\delta_1$ est dirig\'ee par
$(n, n\Chi(V,S) + \Chi(S) - 1)$. D'autre part, pour toute droite rationnelle
$\delta \in [\delta_1,\delta_0\cl[$, on peut trouver un $a \in \op]0,1]$ tel que
$\delta_a$ soit \'egale \`a $\delta$ et que $T_a$ soit un tore convexe d'indice
$2$. En fait, pour cela, il faut \'eventuellement perturber~$\eta$ par une
petite isotopie qui ne d\'etruit pas les propri\'et\'es utiles.

\`A partir de ces observations, la preuve est identique \`a celle de l'assertion
similaire dans la d\'emonstration de la proposition~\ref{p:ghys}.
\end{proof}

\smallskip

Pour terminer la d\'emonstration, on utilise la classification des structures de
contact sur le tore plein \'etablie dans \cite{Gi:bifurcations}. Si $n \Chi(V,S)
= -\Chi(S)$, la droite $\delta_1$ est dirig\'ee par le vecteur $(n,-1)$. Avec
cette condition au bord, le th\'eor\`eme 1.6 de \cite{Gi:bifurcations} affirme
que la restriction de $\xi$ \`a~$W$ est universellement tendue. On se trouve du
coup dans la situation de la section 3.C : $\xi$ est isotope \`a une structure
de contact tangente aux fibres et n'est pas virtuellement vrill\'ee. Si $n=1$ et
si $\Chi(V,S) < -\Chi(S)$, la droite~$\delta_1$ est dirig\'ee par le vecteur
$(1,-m)$, o\`u $m = 1 - \Chi(V,S) - \Chi(S) > 1$. Le th\'eor\`eme 1.6 de~\cite
{Gi:bifurcations} dit alors qu'il y a sur~$W$, \`a isotopie relative au bord
pr\`es, $m-1$~structures de contact tendues qui co\"{\i}ncident avec~$\xi$ sur
$\partial W$, dont une (seule) est universellement tendue. En outre, l'existence
sur~$V$ d'une structure de contact universellement tendue et d'enroulement~$-1$
assure que, si $\xi \res W$ est universellement tendue, $\xi$ l'est aussi. On
obtient ainsi les bornes annonc\'ees.
\end{proof}

\smallskip

\begin{remarque}
La d\'emonstration ci-dessus fait appara\^\i tre que $-1$ est la seule valeur
stric\-tement n\'egative possible pour l'enroulement d'une structure de contact
virtuellement vrill\'ee.

Par ailleurs, les r\'esultats de cette partie permettent de compl\'eter quelque
peu l'\'enonc\'e du th\'eor\`eme~\ref{t:thurston} : si une structure de contact
$\xi$ sur $V$ est d'enroulement positif ou nul, il existe non seulement une
courbe legendrienne isotope \`a la fibre et d'enroulement nul mais tout un tore 
transversal \`a~$\xi$ et dont les caract\'eristiques sont isotopes aux fibres.
Ce tore est la version de contact de la feuille compacte trouv\'ee par W.~{\sc
Thurston}~\cite{Th:these}. Lorsque~$\xi$ est tendue, son existence r\'esulte
soit du th\'eor\`eme~\ref{t:lutz}, soit du lemme~\ref{l:virt}. Lorsque~$\xi$ est
vrill\'ee, un th\'eor\`eme de Y.~{\sc Eliashberg}~\cite{El:overtwist}  assure
qu'une modification de Lutz le long d'une fibre produit une structure de contact
isotope \`a~$\xi$. Or une telle modification fait clairement appara\^\i tre le
tore cherch\'e (voir par exemple~\cite{Gi:bourbaki}).
\end{remarque}

\adresse%
{Emmanuel \textsc{Giroux}\Y
Unit\'e de Math\'ematiques Pures et Appliqu\'ees,\Y
\'Ecole Normale Sup\'erieure de Lyon,\Y
46, all\'ee d'Italie,\Y
69364, Lyon cedex 07, France}

\end{document}